\documentclass[10pt,twocolumn,twoside]{IEEEtran} 
\IEEEoverridecommandlockouts                    
\usepackage{graphicx}
\usepackage{amsmath}
\usepackage{amssymb}
\usepackage{amsfonts}
\usepackage{mathrsfs}
\usepackage{dsfont}
\usepackage{subcaption}
\makeatletter
\def\amsbb{\use@mathgroup \M@U \symAMSb}
\makeatother
\usepackage{color}
\newtheorem{definition}{Definition}
\newtheorem{theorem}{Theorem}
\newtheorem{proposition}{Proposition}

\newtheorem{Remark}{Remark}
\usepackage[usenames,dvipsnames]{xcolor}
\usepackage{epstopdf}
\usepackage[export]{adjustbox}
\DeclareMathOperator{\Ima}{Im}
\newcommand{\norm}[1]{\left\lVert#1\right\rVert}
\newcommand{\vect}[1]{\boldsymbol{#1}} 
\usepackage{multirow}

\usepackage{pifont}
\newcommand{\xmark}{\ding{55}}%

\usepackage{makeidx}
\usepackage{nomencl}
\usepackage{ifthen}
\renewcommand{\nomgroup}[1]{%
\ifthenelse{\equal{#1}{F}}{\item[\textbf{Functions}]}{%
\ifthenelse{\equal{#1}{I}}{\item[\textbf{Functions $\&$ Indices}]}{%
\ifthenelse{\equal{#1}{V}}{\item[\textbf{Variables $\&$ Parameters}]}{%
\ifthenelse{\equal{#1}{S}}{\item[\textbf{Sets}]}{}}}}
}
\makenomenclature

\renewcommand{\nompreamble}{Main symbols used within the paper.
}

\allowdisplaybreaks

\usepackage{balance}
\newtheorem{lemma}{Lemma}
\newtheorem{problem}{Problem}
\newtheorem{design}{Design Condition}
\usepackage{enumerate}

\title{Privacy of distributed optimality schemes in power networks}
\author{{{Andreas Kasis},\thanks{{Andreas Kasis and Kanwal Khan are first authors who contributed equally.}} {Kanwal Khan},\thanks{Andreas Kasis, Kanwal Khan,  Marios M. Polycarpou and Stelios Timotheou are with the KIOS Research and Innovation Center of Excellence and the Department of Electrical and Computer Engineering, University of Cyprus, Cyprus; e-mails:  kasis.andreas@ucy.ac.cy, {kanwalhasan26@gmail.com}, mpolycar@ucy.ac.cy, timotheou.stelios@ucy.ac.cy.} Marios M. Polycarpou and Stelios Timotheou \thanks{This work was funded by the European Union’s Horizon 2020 research and innovation program under grant agreements No. 891101 (SmarTher Grid) and No. 739551 (KIOS CoE), and from the Republic of Cyprus through the Directorate General for European Programs, Coordination, and Development.}
\thanks{A preliminary version of this work will appear in \cite{Privacy_CDC}. 
This work extends the scheme presented in \cite{Privacy_CDC}, introduces a new scheme,  and includes additional discussion and analytic and simulation results.
}
}}
\begin{document}
\date{}
\maketitle

\begin{abstract}
The increasing participation of local generation and controllable demand units within the power network motivates the use of distributed schemes for their control. Simultaneously, it raises two issues;  achieving an optimal power allocation among these units, and securing the privacy of the generation/demand profiles. This study considers the problem of designing distributed optimality schemes that preserve the privacy of the generation and controllable demand units within the secondary frequency control timeframe. We propose a consensus scheme that includes the generation/demand profiles within its dynamics, keeping this information private when knowledge of its internal dynamics is not available. However, the prosumption profiles may be inferred using knowledge of its internal model. We resolve this by proposing a privacy-preserving scheme which ensures that the generation/demand 
cannot be inferred from the communicated signals. For both proposed schemes, we provide analytic stability, optimality and privacy guarantees and show that the secondary frequency control objectives are satisfied. The presented schemes are distributed, locally verifiable and applicable to arbitrary network topologies. Our analytic results are verified with simulations on a 140-bus system, where we demonstrate that the proposed schemes offer enhanced privacy properties, enable an optimal power allocation and preserve the stability of the power network.
\end{abstract}

\nomenclature[V]{{$\omega_j$}}{ {frequency {deviation at bus $j$}}}
\nomenclature[V]{{$\eta_{ij}$}}{ {power angle difference between bus $i$ and bus $j$}}
\nomenclature[V]{{$p^M_{k,j}$}}{ {$k$th mechanical power injection {at bus $j$}}}
\nomenclature[V]{{$x_{k,j}$}}{ {internal states of $k$th mechanical power injection unit {at bus $j$}}}
\nomenclature[V]{{$d^c_{k,j}$}}{ {demand of the $k$th controllable load {at bus $j$}}}
\nomenclature[V]{{$p^L_{k,j}$}}{ {$k$th uncontrollable demand unit {at bus $j$}}}
\nomenclature[V]{{$p_{ij}$}}{ {power transfer from bus $i$ to bus $j$} }
\nomenclature[V]{{$B_{ij}$}}{ {line susceptance {between bus $i$ and bus $j$}}}
\nomenclature[V]{{$M_j$}}{ {generator inertia {at bus $j$}}}
\nomenclature[V]{{$D_j$}}{ {frequency damping {at bus $j$}}}
\nomenclature[V]{{$p^c_j$}}{ {power command {at bus $j$}}}
\nomenclature[V]{{$p^c_{k,j}$}}{ {power command variable {at $k$th prosumnption unit at bus $j$}}}
\nomenclature[V]{{$\psi_{i,j}$}}{ {integral of power command difference between units $i$ and $j$}}
\nomenclature[V]{{$\hat{H}, H$}}{ {incidence matrices associated with the communication graphs of the \textit{Primal-Dual}, and the \textit{Extended Primal-Dual} and \textit{Privacy-Preserving} schemes respectively}}
\nomenclature[S]{{$\mathcal{N}$}}{ {set of buses}}
\nomenclature[S]{{$\mathcal{E}$}}{ {set of transmission lines}}
\nomenclature[S]{{$\mathcal{N}^G_j$}}{ {set of generation units at bus $j$}}
\nomenclature[S]{{$\mathcal{N}^L_j$}}{ {set of controllable demand units at bus $j$}}
\nomenclature[S]{{$\mathcal{N}_j$}}{ {set of uncontrollable demand units at bus $j$}}
\nomenclature[S]{{$\mathcal{N}^p_{j}$}}{ {set of buses preceding bus $j$ in the power network}}
\nomenclature[S]{{$\mathcal{N}^s_{j}$}}{ {set of buses succeeding bus $j$ in the power network}}
\nomenclature[I]{{$x^*$}}{ {equilibrium point of variable $x$}}
\nomenclature[I]{{$\dot{x}$}}{ {time derivative of function of time $x$}}
\nomenclature[S]{{$\hat{\mathcal{E}}$}}{ {set of communication lines in the \textit{Primal-Dual} scheme}}
\nomenclature[S]{{$\widetilde{\mathcal{E}}$}}{ {set of communication lines in the \textit{Extended Primal-Dual} and \textit{Privacy-Preserving} schemes}}
\nomenclature[S]{{$\widetilde{\mathcal{N}}$}}{ {set of prosumption units}}
\nomenclature[V]{{$\widetilde{s}_j$}}{ {vector with all generation and controllable and uncontrollable demand units {at bus $j$}}}
\nomenclature[V]{{$n$}}{ {privacy-enhancing signal}}
\nomenclature[V]{{$\zeta_j$}}{ {aggregate demand minus generation at bus $j$}}
\printnomenclature

\section{Introduction}

\textbf{Motivation and literature survey:}
The increasing penetration of renewable sources of generation is expected to cause more frequent generation-demand imbalances within the power network, which may harm power quality and even cause blackouts \cite{ipakchi2009grid}.
Controllable demand is considered  to be a means to address this issue, since loads may provide a fast response {to counterbalance} intermittent generation \cite{kamyab2015demand}.
However,  the increasing  number of such active units
makes traditionally implemented centralized {control} schemes expensive and inefficient, motivating the adoption of distributed schemes.
Such schemes offer many advantages, such as scalability, reduced expenses associated with the necessary communication infrastructure {and enhanced reliability due to the absence of a single point of failure.}

The introduction of  controllable loads and local renewable generation
 raises an issue {of economic optimality} in the  power allocation.
In addition, the introduction of smart meters for the monitoring of {generation and demand} units poses a privacy threat for the citizens, since readings may be used to expose customers daily life and habits, by inferring the users energy consumption patterns and types of appliances \cite{zeifman2011nonintrusive}. 
{For example, this issue} led the Dutch Parliament  to prohibit the deployment of smart meters until the privacy concerns are resolved \cite{erkin2013privacy}, {as well as} several counties and cities in California to vote for making smart meters illegal in their jurisdictions \cite{hess2014wireless}.
These {concerns} motivate the design of distributed schemes that will simultaneously achieve an optimal power allocation and preserve the privacy of local prosumption profiles.

In recent years, various studies considered the use of  decentralized/distributed control schemes for generation and controllable demand with applications to both primary \cite{devane2016primary}, \cite{kasis2017stability}, \cite{zhao2014design}, \cite{kasis2021primary} and secondary  \cite{kasis2016primary}, \cite{trip2016internal}, \cite{li2015connecting}, \cite{chen2020distributed} frequency regulation, where the objectives are to ensure generation-demand balance and that the frequency attains its nominal value at steady state respectively.
In addition, the problem of obtaining an optimal power allocation within the secondary frequency control timeframe has received broad  attention in the literature \cite{mallada2017optimal}, \cite{zhao2015distributed}, \cite{zhao2018distributed}.
These studies considered suitably constructed optimization problems and designed the system equilibria to coincide with the solutions to these problems.
In many studies, the control dynamics were inspired from the dual of the considered optimization problems   \cite{li2015connecting}, \cite{low2014distributed},  \cite{kasis_TCST}.
Such schemes,  usually referred to in the literature as \textit{Primal-Dual schemes},  yield an optimal power allocation and at the same time {enable the satisfaction of} operational constraints.
Alternative distributed schemes, which ensure that frequency attains its nominal value at steady state by using the generation outputs, have also been proposed \cite{kasis2020distributed}, \cite{trip2017distributed}.
However, the use of real-time knowledge of the generation and controllable demand in the proposed schemes  may compromise the privacy of prosumers.

The topic of preserving the privacy of generation and demand units {has recently} attracted wide attention in the literature.
Different types of privacy concerns, resulting from the integration of information and communication technologies in the smart grid, are mentioned in \cite{zeadally2013towards}. 
In addition, \cite{siddiqui2012smart} analyzes various smart grid privacy issues and discusses recently proposed solutions for enhanced privacy, while \cite{yu2013privacy}  proposes a privacy-preserving power request scheme. 
In addition, \cite{fioretto2019differential} uses the differential privacy framework  to provide privacy guarantees and \cite{eibl2017differential} studies the effect of differential privacy on smart metering data.
Moreover,  homomorphic encryption has been used in  \cite{marmol2012not} to enable the direct connection and exchange of data between electricity suppliers and final users, while preserving the privacy in the smart grid. 
A privacy-preserving aggregation scheme is proposed in \cite{lu2012eppa} which considers various security threats.
{A network equivalent approach is developed in  \cite{zheng2021dynamic} which preserves the information privacy of  integrated electricity and heat systems.} 
 The use of energy storage units to preserve the privacy of user consumption has been considered in \cite{yang2014cost} and \cite{zhang2016cost}.
 Furthermore, \cite{wu2021privacy} and \cite{dvorkin2020differentially} aim to simultaneously preserve the privacy of individual agents and enable an optimal power allocation  using homomorphic encryption and differential privacy respectively.
 Both approaches result in suboptimal allocations, which suggests a trade-off between optimality and privacy. 
Several existing techniques {that aim at} preventing disclosure of private data are also discussed in  \cite{souri2014smart}.

Although the problems of preserving the privacy of power prosumption and obtaining an optimal power allocation in power networks have been independently studied, 
{to the authors best knowledge, no study has managed to simultaneously achieve these objectives.}
{In addition, the impact of such schemes  on the stability and {dynamic} performance of the power grid has not been investigated.}
This study aims to jointly consider these objectives within the secondary frequency control timeframe.

\textbf{Contribution:}
This paper  studies the problem of  providing optimal frequency regulation within the secondary frequency control timeframe while preserving the privacy of generation and controllable demand profiles.
We first propose an optimization problem that ensures that secondary frequency regulation objectives, i.e. achieving generation-demand balance and frequency attaining its nominal value at steady state, are satisfied.
In addition, to facilitate the interpretation of our privacy results, we define two types of eavesdroppers;
(i) \textit{naive eavesdroppers}, that do not {possess/make use of knowledge of the system dynamics to analyze} the intercepted information and (ii) \textit{informed or intelligent eavesdroppers} that {use} knowledge of the underlying system dynamics to infer the prosumption profiles.

We consider a distributed scheme that has been extensively studied in the literature, usually referred to as the \textit{Primal-Dual scheme}, that enables an optimal power allocation and the satisfaction of system constraints, and explain why it {causes} privacy issues.
Inspired by the \textit{Primal-Dual scheme}, we propose the \textit{Extended Primal-Dual scheme} that incorporates a distributed controller at each privacy-seeking unit of the power grid.
The latter replaces the communication of prosumption profiles with a consensus signal providing  privacy  against naive eavesdroppers.
However, we explain how intelligent eavesdroppers may infer the prosumption profiles using the communicated signal trajectories and knowledge of the underlying system dynamics.
To resolve this, we propose the \textit{Privacy-Preserving scheme}, which 
incorporates two important features into the \textit{Extended Primal-Dual scheme}, such that privacy against intelligent eavesdroppers is achieved.
In particular, the proposed scheme {continuously} alters the speed of response of each controller, making model based inference inaccurate. Moreover, it adds bounded noise to the prosumption information within each controller, with a maximum magnitude proportional to the local frequency deviation. 
The latter yields changes in all controllers when a disturbance occurs, making it hard to detect the origin of the disturbance.
These properties
ensure that  the \textit{Privacy-Preserving scheme} {guarantees} the  privacy of the prosumption units against intelligent eavesdroppers.
{On the other hand, due to its additional features,  the \textit{Privacy-Preserving scheme}
could  potentially result in  slower convergence, since the controllers response speed is reduced.}
For both proposed schemes, we provide analytic stability guarantees and show that an optimal power allocation is achieved at steady state.
In addition, the proposed schemes are distributed and applicable to arbitrary network topologies, {while} the proposed conditions are locally verifiable. 

Our analytic results are {illustrated} with numerical simulations on the NPCC 140-bus system which validate that the proposed schemes enable an optimal power allocation and satisfy the secondary frequency regulation objectives.
In addition, we demonstrate how the \textit{Extended Primal-Dual} and the \textit{Privacy-Preserving} schemes offer privacy of the prosumption profiles  against naive and intelligent eavesdroppers respectively.

To the authors best knowledge, this is the first study that:
\begin{enumerate}[(i)]
\item Jointly studies the {privacy, optimality and stability}  properties of distributed schemes within the secondary frequency control timeframe.
\item Proposes distributed schemes that yield an optimal power allocation  and simultaneously preserve the privacy  of the {prosumption} profiles.
In particular, the proposed schemes offer privacy guarantees against naive (\textit{Extended Primal-Dual scheme}) and informed (\textit{Privacy-Preserving scheme})  eavesdroppers respectively. For the proposed schemes, we show that stability is guaranteed and that the secondary frequency control objectives are satisfied.
\end{enumerate}

\textbf{Paper structure:}
In Section \ref{II} we present the dynamics of the power network,  the considered optimization problem and the problem statement. 
{In Section \ref{Sec_Primal_Dual} we examine the \textit{Primal-Dual} scheme and discuss its privacy issues.}
In Sections \ref{III} and \ref{IV} we present the proposed \textit{Extended Primal-Dual} and \textit{Privacy-Preserving} schemes respectively  and provide our main analytic results.
In Section \ref{V} we  validate our main results through numerical simulations on the NPCC 140-bus system. 
Finally, conclusions are drawn in Section \ref{VI}.
The proofs of the main analytic results {(Theorems \ref{theorem1}, \ref{theorem2}, Propositions \ref{prop_privacy_naive}--\ref{proposition2} and Lemma \ref{Lemma1}) are omitted due to space restrictions and are provided in \cite{arXiv_privacy}}.

\begin{figure}
\includegraphics[trim = 10mm 7mm 0 3mm, clip = true, scale=0.95]{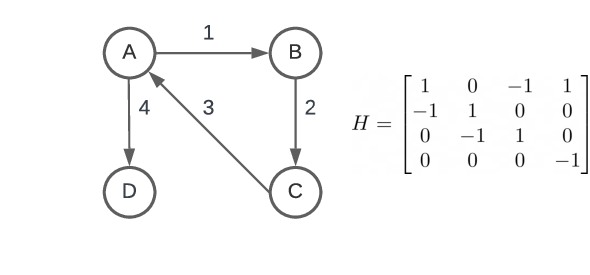}
\vspace{-2mm}
\caption{{The incidence matrix $H$ of a simple $4$-node graph.}}
\label{Fig_incidence_matrix}
\vspace{-4mm}
\end{figure}

\textbf{Notation:}
Real numbers and the set of n-dimensional vectors with real entries are denoted by $\mathbb{R}$  and $\mathbb{R}^n$ respectively.
The $p$-norm of a vector $x \in \mathbb{R}^n$ is given by $\norm{x}_p = (|x_1|^p + \dots + |x_n|^p)^{1/p}, 1 \leq p < \infty$.
A function $f: \mathbb{R}^n \rightarrow \mathbb{R}^m$ is said to be locally Lipschitz continuous at $x$ if there exists some neighbourhood $X$ of $x$ and some constant $L$ such that $\norm{f(x) - f(y)} \leq L \norm{x - y}$ for all $y \in X$, where $\norm{.}$ denotes any $p$-norm.
A matrix $A \in \mathbb{R}^{n \times n}$ is called diagonal if $A_{ij} = 0$ for all $i \neq j$. In addition, $A \preceq 0$ indicates that the matrix $A$ is  negative semi-definite. 
The image of a vector $x$ is denoted by $\Ima(x)$.
The cardinality of a discrete set $\mathcal{S}$ is denoted by $|\mathcal{S}|$.
{A set $\mathcal{B}$ is a proper subset of a set $\mathcal{A}$ if $\mathcal{B}  \subset \mathcal{A}$ and $\mathcal{B}  \neq \mathcal{A}$.}
For a graph with sets of nodes and edges denoted by  $\mathcal{A}$ and $\mathcal{B}$ respectively, 
we define the incidence matrix $H  \in \mathbb{R}^{|{\mathcal{A}}| \times |{\mathcal{B}}|}$ as follows
\begin{gather*}
H_{ij} = \begin{cases}
+1, \text{ if } i \text{ is the positive end of edge } j \in \mathcal{B}, \\
-1, \text{ if } i \text{ is the negative end of edge } j \in \mathcal{B}, \\
0, \text{ otherwise.}
\end{cases}
\end{gather*}
{An illustrative example of the incidence matrix of a simple graph is presented in Fig. \ref{Fig_incidence_matrix}.}
We use $\vect{0}_n$ and $\vect{1}_n$  to denote $n$-dimensional vectors with all elements equal to $0$ and $1$ respectively.
Finally, for a state $x \in \mathbb{R}^n$, we let $x^*$ denote its equilibrium value.

\section{Problem Formulation}\label{II}

\subsection{Power network model}\label{Sec_power_network}
 We describe the power network by a connected graph ${\mathcal{(N,E)}}$ where ${\mathcal{N}}=\{1,2,..,|{\mathcal{N}}|\}$ is the set of buses and ${\mathcal{E}}\subseteq {\mathcal{N}}\times {\mathcal{N}}$ the set of transmission lines connecting the buses. The term $(i,j)$ denotes the link connecting buses $i$ and $j$. The graph ${\mathcal{(N,E)}}$ is assumed to be directed with an arbitrary direction, so that if $(i,j)\in{{\mathcal{E}}}$ then $(j,i)\notin {{\mathcal{E}}}$. 
 For each $j \in \mathcal{N}$, we define the sets of predecessor and successor buses by  $\mathcal{N}^p_{j} = \{k : (k,j) \in \mathcal{E}\}$ and $\mathcal{N}^s_{j} = \{k : (j,k) \in \mathcal{E}\}$ respectively. 
It should be noted that the form of the considered dynamics is unaffected by changes in the graph ordering and the results presented in this paper are independent of the choice of direction.
The following assumptions are made for the network: \newline
1) Bus voltage magnitudes are $|V_j| = 1$ per unit for all $j \in \mathcal{N}$. \newline
2) Lines $(i,j) \in \mathcal{E}$ are lossless and characterized by the magnitudes of their susceptances $B_{ij} = B_{ji} > 0$. \newline
3) Reactive power flows do not affect bus voltage phase angles and frequencies. \newline
4) The relative phase angles are sufficiently small such that the approximation $\sin \eta_{ij} = \eta_{ij}$ is valid. \newline
The first three assumptions  {are standard in secondary frequency regulation studies}  \cite{kasis2017stability, zhao2014design, kasis2016primary, li2015connecting}, {\cite{chen2020distributed, yang2021decentralized, baros2021examining, wang2022frequency}}. 
{These assumptions usually hold in medium to high voltages, associated with transmission systems, where lines are dominantly inductive and voltage variations are small and tightly controlled.}
The fourth assumption is valid when the network operates in nominal conditions, where relative phase angles are small\footnote{It should be noted that the results presented in this paper can be extended by considering sinusoidal phase angles, see e.g. the approach in \cite{kasis2017stability}. We have opted not to consider this case for simplicity and to keep the main focus of the paper on the privacy aspects of the proposed schemes.}.
It should be noted that the theoretical results presented in this paper are validated with numerical simulations in Section \ref{V}, on a comprehensive power network model.

We use the swing equations to describe the rate of change of frequency at buses \cite{machowski2020power}. 
In particular, at each bus we consider a set of generation and controllable and uncontrollable demand units.
This motivates the following system dynamics: 
\begin{subequations}\label{eq1}
\begin{align}
    \dot{\eta}_{ij} &=\omega_{i}-\omega_{j}, (i,j)\in{\mathcal{E}}, \label{eq1a} \\
  \begin{split}
M_{j}\dot{\omega}_{j} & = \sum_{k \in \mathcal{N}^G_j} p^M_{k,j} - \sum_{k \in \mathcal{N}^L_j}  d^c_{k,j} - \sum_{k \in \mathcal{N}_j} p^L_{k,j}  - D_{j}\omega_{j}\\
& - \sum_{i \in \mathcal{N}^s_j} p_{ji} +
   \sum_{i \in \mathcal{N}^p_j} p_{ij}, j\in{\mathcal{N}}, \label{eq1b}
\end{split} \\
 p_{ij} &=B_{ij}\eta_{ij},(i,j)\in{\mathcal{E}}. \label{eq1c}
\end{align}
\end{subequations}
In system \eqref{eq1}, variable $\omega_{j} $ represents the deviation of the frequency at bus $j$ from its nominal value, namely 50 Hz (or 60 Hz).
Variable 
$p^M_{k,j}$ represents the mechanical power injection associated with the $k$th generation unit at bus $j$.
Moreover, $d^c_{k,j}$ denotes the demand associated with the $k$th controllable load at bus $j$. $\mathcal{N}^G_j$ and $\mathcal{N}^L_j$ represent the sets of generation units and controllable loads, which are jointly referred to as active elements or active units, at bus $j$ respectively.
Each of these units are associated with a privacy-seeking user or entity.
The set of active units  at bus $j$ is given by $\mathcal{N}_j = \mathcal{N}^G_j \cup \mathcal{N}^L_j$.
The variable $p^L_{k,j}$ represents the uncontrollable demand associated with the $k$th active unit at bus $j$. 
 Furthermore, the time-dependent variables $\eta_{ij} $ and $p_{ij}$   represent, respectively, the power angle difference and the power transmitted from bus ${i}$ to bus ${j}$. 
 The quantities $B_{ij}$  represent the line susceptances between buses $i$ and $j$.
 Finally, the positive constants $D_{j}$ and $M_j$ represent the generation damping and inertia at bus $j$ respectively. 
 The generation and consumption will be jointly referred to as prosumption.

\begin{Remark}
An alternative, but equivalent, representation of \eqref{eq1} could include 
a single variable at each bus representing the aggregation of uncontrollable demand.
We opted to associate uncontrollable loads with active units   to facilitate the study of their privacy properties. 
The benefits of this representation are evident in Sections \ref{III} and \ref{IV}.
Note that when no uncontrollable load is associated with some generation or controllable demand unit, then $p^L_{k,j} = 0$. 
\end{Remark}

\subsection{Generation and controllable demand dynamics}

 We will study the behavior of the power system under the following dynamics for generation and controllable loads,
\begin{subequations}\label{eq2}
\begin{align}
\tau_{k,j} \dot{x}_{k,j} &= -x_{k,j} + m_{k,j}(u_{k,j} - \omega_j), k \in \mathcal{N}^G_j, j\in{\mathcal{N}}, \label{eq2a} 
\\
 p^M_{k,j} &= x_{k,j} + h_{k,j} (u_{k,j} - \omega_j), k \in \mathcal{N}^G_j, j\in{\mathcal{N}}, \label{eq2b} \\ 
 d^c_{k,j} &= - h_{k,j} (u_{k,j} - \omega_j), k \in \mathcal{N}^L_j, j \in{\mathcal{N}}, \label{eq2c}
               \end{align}             
  \end{subequations}
where $x_{k,j} \in \mathbb{R}$ represents the internal state, and $\tau_{k,j} > 0$ and $m_{k,j} > 0$ the time and droop constants associated with generation unit $k$ at bus $j$ respectively.
The positive constant $h_{k,j}$ represents the damping associated with active unit $k$ (generation or controllable load) at bus $j$.
In addition, $u_{k,j}$ represents the control input to the $k$th active unit at bus $j$, the dynamics of which are discussed in the following sections.
{It should be noted that generation and controllable demand units, as well as their inputs, evolve in continuous time.}

We consider first-order generation dynamics and static controllable demand for simplicity and to keep the focus of the paper on developing a privacy-preserving scheme. 
More involved generation and demand dynamics  could be considered by applying existing results (e.g. \cite{kasis2017stability},  \cite{kasis2016primary}, \cite{kasis_TCST}).

For convenience, we define the vectors $p^M_j = [p^M_{k,j}]_{k \in \mathcal{N}^G_j}$,  $d^c_j = [d^c_{k,j}]_{k \in \mathcal{N}^L_j}$, $p^L_j = [p^L_{k,j}]_{k \in \mathcal{N}_j}$, $p^M = [p^M_j]_{j \in \mathcal{N}}$, $d^c = [d^c_j]_{j \in \mathcal{N}}$ and $p^L = [p^L_j]_{j \in \mathcal{N}}$.

\subsection{{Prosumption cost minimization} problem}\label{sec_optimization}

In this section we form an optimization problem that aims to minimize the costs associated with generation and controllable demand and simultaneously achieve generation-demand balance.
The considered optimization problem is described below.

A cost $\frac{1}{2}q_{k,j}(p^M_{k,j})^2$ is incurred when the generation unit $k$ at bus $j$   produces a power output of $p^M_{k,j}$. In addition, a cost $\frac{1}{2}q_{k,j}(d^c_{k,j})^2$ is incurred when controllable load $k$ at bus $j$ adjusts its demand to $d^c_{k,j}$.
{It should be noted that quadratic cost functions are considered since those can locally approximate general convex cost functions.}
The optimization problem is to obtain the vectors $p^M$ and $d^c$ that minimize the  cost associated with the aggregate generation and controllable demand and simultaneously achieve power balance.
The considered optimization problem is presented below.
\begin{equation}\label{eq3}
\begin{aligned}
\min_{p^M,d^c}  & \sum_{j \in \mathcal{N}} ( \sum_{k \in \mathcal{N}^G_j} \frac{1}{2}q_{k,j}(p^M_{k,j})^2 + \sum_{k \in \mathcal{N}^L_j} \frac{1}{2}q_{k,j}(d^c_{k,j})^2)\\
\text{ subject to }
&  \sum_{j \in \mathcal N} (\sum_{k \in \mathcal{N}^G_j} p^M_{k,j} - \sum_{k \in \mathcal{N}^L_j} d^c_{k,j} - \sum_{k \in \mathcal{N}_j} p^L_{k,j}) = 0.
 \end{aligned}
 \end{equation}
The equality constraint in \eqref{eq3} requires all the uncontrollable loads to be matched by the generation and controllable demand,  such that generation-demand balance is achieved.
The equality constraint also guarantees that  the frequency attains its nominal value at equilibrium, which is a main objective of secondary frequency control.
The latter follows  by summing \eqref{eq1b} at steady state over all buses, which yields $\sum_{j \in \mathcal{N}} D_j \omega_j = 0$, and noting that  frequency synchronizes at equilibrium from \eqref{eq1a}.

{
\begin{Remark}\label{rem_constraints}
The optimization problem \eqref{eq3} does not consider any power line constraints.
The incorporation of such constraints has been studied in the literature \cite{chen2020distributed, zhao2015distributed}, where suitable approaches have been developed.
We opted not to include power line constraints in this study to keep its focus on the privacy aspect of the control design.
\end{Remark}
}

\subsection{{Eavesdropper and privacy definitions}}\label{Sec_eavesdroppers}

{In this section, we define the two considered eavesdropper types, inspired from \cite{parsaeefard2015improving},  and present two notions of privacy to facilitate the interpretation and intuition of our results.}

\begin{definition}\label{Def_eavesdropper}
{An eavesdropper is a person or entity that aims to extract private information  by intercepting the signals communicated to and from generation and controllable demand units.}
Eavesdroppers are classified  as follows:
\begin{enumerate}[(i)]
    \item Naive eavesdroppers, who 
    {posses knowledge of: \newline 
    (K1) All signals communicated to and from a given unit, for which it aims to obtain private information.}
    \item {Informed or intelligent eavesdroppers, who posses knowledge of K1 and: \newline
    (K2) The underlying control dynamics of the system.}
\end{enumerate}
\end{definition}

Definition \ref{Def_eavesdropper} presents two types of eavesdroppers, based on {whether they make use of knowledge of the underlying system dynamics to infer private information.}
In particular, naive eavesdroppers 
{have no knowledge of the system model that  {may be utilized to} analyze the} intercepted signals. 
They only try to overhear sensitive information. 
Informed eavesdroppers  {analyze the intercepted signals using knowledge of} the underlying dynamics. It is intuitive to note that privacy against intelligent eavesdroppers implies privacy against naive eavesdroppers but not vice versa.

{
\begin{Remark}\label{remark_eavesdroppers}
It should be noted that naive and intelligent eavesdroppers are assumed to  possess knowledge of K1 and K1 and K2 respectively throughout the considered time duration. 
In addition, it is assumed that naive and intelligent eavesdroppers stay as such throughout the considered time duration, i.e. knowledge of $K2$ is not gained or lost by eavesdroppers during the considered timeframe.
\end{Remark}
}

{
Below we provide a definition of a private prosumption trajectory and profile, used throughout the rest of the manuscript. We remind that $s^*$ denotes the equilibrium value of $s$, i.e. $s^* = \lim_{t \rightarrow \infty} s(t)$.

\begin{definition}\label{dfn_privacy}
The following two notions of prosumption privacy are considered: \newline 
(i) A prosumption trajectory is called private against an eavesdropper type if the knowledge available to the eavesdropper does not allow the estimation of $s(t), t \geq 0, s \neq s^*$.  \newline
(ii) A prosumption profile is called private against an eavesdropper type if the knowledge available to the eavesdropper does not allow the estimation of its trajectory and steady state values, i.e. of $s(t), t \geq 0$.
\end{definition}

The considered privacy definition implies that a prosumption trajectory is private when  an eavesdropper cannot accurately estimate its initial condition and values when not at steady state. 
The privacy of a prosumption profile requires in addition the privacy of its steady state value. 
The distinction between the two notions {enables} privacy guarantees based on different conditions (see Section  \ref{sec_privacy_analysis}).
}

\subsection{Problem Statement}

This paper aims to design control schemes that enable  stability and optimality guarantees and at the same time preserve the privacy of all active units.
The problem is stated below.

\begin{problem}\label{problem_definition}
Design a control scheme that:
\begin{enumerate}[(i)]
\item Preserves the privacy {of the prosumption profiles}  against intelligent eavesdroppers.
\item Enables asymptotic stability guarantees.
\item Uses local information and locally verifiable conditions.
\item Yields an optimal steady-state power allocation.
\item Applies to arbitrary connected network configurations.
\end{enumerate}
\end{problem}

Problem \ref{problem_definition} aims 
to design a control scheme  that enables stability guarantees, ensures an optimal power allocation at steady state,
and guarantees the privacy of the generation/demand profiles against informed eavesdroppers, {following Definitions \ref{Def_eavesdropper} and \ref{dfn_privacy}.}
In addition, we aim to design a scheme that relies on locally available information and locally verifiable conditions, to  enable scalable designs. Finally, it is desired that the proposed scheme is applicable to general network topologies.

\section{Primal-Dual scheme}\label{Sec_Primal_Dual}

{In this section we
examine} a distributed scheme that has been widely studied in the literature \cite{kasis2017stability}, \cite{li2015connecting}, \cite{zhao2015distributed}, \cite{low2014distributed},     usually referred to as the \textit{Primal-Dual scheme},  that  enables an optimal power allocation, and discuss its resulting privacy issues.

To describe the \textit{Primal-Dual scheme}, we consider a connected communication graph ($\mathcal{N},\hat{\mathcal{E}}$), where $ \hat{\mathcal{E}}$ represents the set of communication lines among the buses, i.e. $(i,j) \in {\hat{\mathcal{E}}}$  if buses $i$ and $j$ communicate. 
In addition, we let $\hat{H}$ be the incidence matrix of  ($\mathcal{N},\hat{\mathcal{E}}$) and define the variable $\zeta_j = \vect{1}^T_{|\mathcal{N}_j|}p^L_{j}  + \vect{1}^T_{|\mathcal{N}^L_j|}d^c_{j}  -\vect{1}^T_{|\mathcal{N}^G_j|}p^M_{j}$ for all $j \in \mathcal{N}$.
The prosumption input dynamics are given by
\begin{subequations}\label{eq4}
\begin{align}
\hat{\Gamma} \dot{\psi} &= \hat{H}^T p^c, \label{eq4a}
    \\
    \bar{\Gamma} \dot{p}^c &= \zeta - \hat{H} \psi,\label{eq4b}
\\
  u_{k,j} &= p^c_j, k \in \mathcal{N}_j, j \in{\mathcal{N}},\label{eq4c}
\end{align}
\end{subequations}
where the diagonal matrices $\hat{\Gamma} \in \mathbb{R}^{|\hat{\mathcal{E}}| \times |{\hat{\mathcal{E}}|}}$ and  $\bar{\Gamma} \in \mathbb{R}^{|{\mathcal{N}}| \times |{\mathcal{N}}|}$ contain the positive time constants associated with \eqref{eq4a} and \eqref{eq4b} respectively
and $p^c_j$ is a power command variable associated with bus $j$ and shared with communicating buses. In addition, variable  $\psi$ is a state of the \textit{Primal-Dual scheme} that integrates the difference in power command variables between communicating buses.
The input for all active elements at bus $j$ is given by the local power command value $p^c_j$, via \eqref{eq4c}.

{The Primal-Dual scheme aims to achieve two objectives through its dynamics.
Firstly, it aims to ensure} that the secondary frequency control objectives, i.e. ensuring generation/demand balance and the frequency attaining its nominal value, are satisfied at steady state.
{The latter follows by summing \eqref{eq1b} and \eqref{eq4b} at steady state over all $j \in \mathcal{N}$, which yield $\sum_{j \in \mathcal{N}} (p^L_j + d^{c,*}_j + D_j \omega^*_j - p^{M,*}_j) = 0$ and $\sum_{j \in \mathcal{N}} \zeta^*_j = \sum_{j \in \mathcal{N}} (p^L_j + d^{c,*}_j - p^{M,*}_j) = 0$ respectively, which}
 implies that  $\omega^* = \vect{0}_{|\mathcal{N}|}$ from the  synchronization of frequency at equilibrium, as follows from \eqref{eq1a}. 
The dynamics in \eqref{eq4a} enable the synchronization of the power command variables at steady state. This property is useful to provide an optimality interpretation of the system's equilibria.
It should be noted that the stability and optimality of the \textit{Primal-Dual scheme} \eqref{eq4} for a wide class of generation/demand dynamics, including those in \eqref{eq2}, have been analytically shown in the literature (e.g. \cite{kasis2017stability}).

\begin{Remark}\label{remark_privacy}
A shortcoming of the \textit{Primal-Dual scheme} \eqref{eq4} is the requirement for real-time knowledge of the generation and demand from all active units in the network.
In practice, this would require the transmission of this information to a central controller at each bus, in order to calculate $\zeta_j$, exposing the local generation/demand profiles 
to a naive eavesdropper who intercepts these signals.
The latter compromises the  privacy of the prosumption profiles.
\end{Remark}

\section{Extended Primal-dual scheme}\label{III}

{In this section, we aim to resolve the privacy issues associated with the \textit{Primal-Dual scheme}.
To achieve this we propose the \textit{Extended Primal-Dual scheme}, which enables the privacy  of the prosumption profiles  against naive eavesdroppers
 and simultaneously retains the beneficial properties of the \textit{Primal-Dual scheme} in terms of being distributed and achieving an optimal power allocation at steady state.}

\subsection{Extended Primal-Dual scheme}\label{sec_extended_PD}

\begin{figure*}[ht!]
\centering\includegraphics[width=\textwidth,height=\textheight,keepaspectratio]{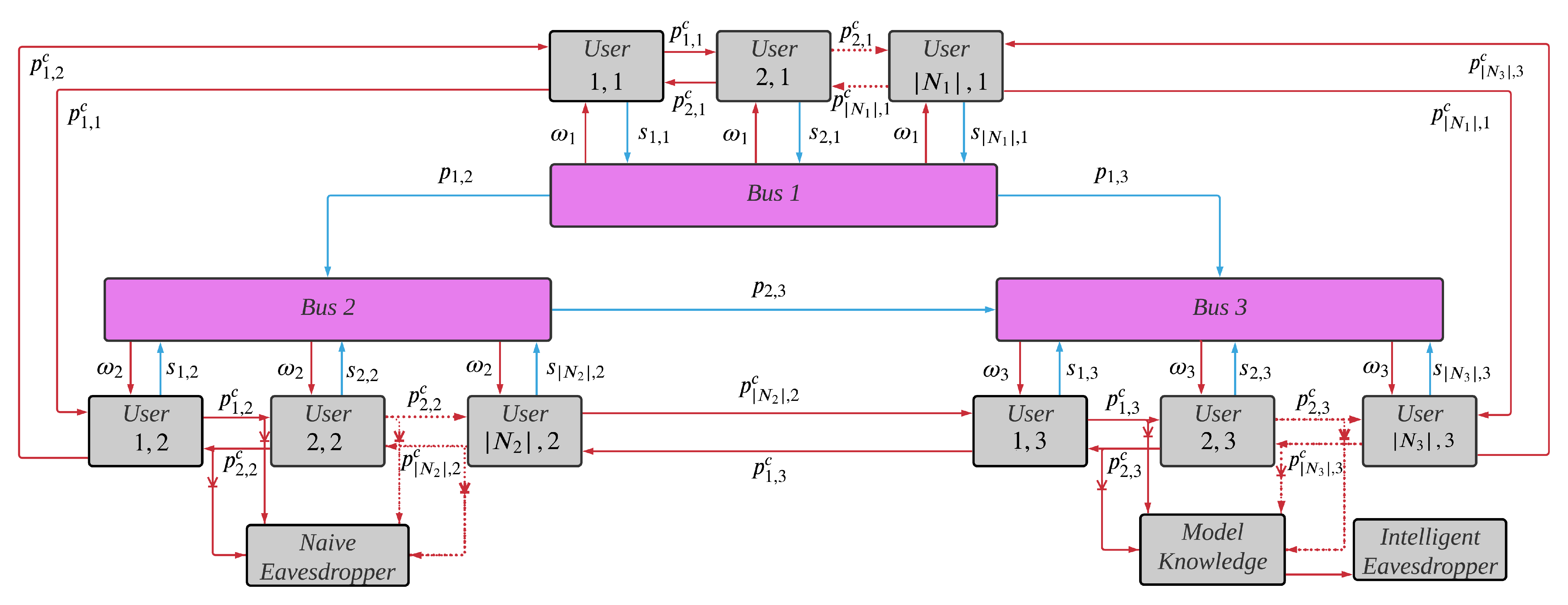}
  \vspace{-4mm}
  \caption{Schematic representation of  system \eqref{eq1}, \eqref{eq2}, \eqref{eq_s}, \eqref{eq5} on a simple 3-bus network. Privacy-seeking users are associated with prosumption units. Blue lines represent  power transfers whereas red lines represent information flows.
  {Red dotted lines represent the communication between all users which can either be direct or indirect through a connected communication network.}
  {Users monitor the local frequency and communicate their respective power command values to neighbouring users.}
  Both naive and intelligent eavesdroppers intercept the communicated signals between users, but only intelligent eavesdroppers 
  {possess knowledge of the underlying dynamics, that may be used to analyze the intercepted information.}}
  \label{Fig1}
  \vspace{-4mm}
\end{figure*}

In this section we present a scheme that aims to improve the privacy properties of the generation/demand profiles. 
In contrast to \eqref{eq4}, which includes a controller at each bus, the {proposed} scheme employs a controller at each privacy-seeking unit (generator or controllable load).
We demonstrate that the presented scheme offers  privacy  against naive  eavesdroppers and simultaneously enables an optimal power allocation.

To describe the new scheme, we consider a communication network characterized by a connected graph ($\widetilde{\mathcal{N}}, \widetilde{\mathcal{E}}$), where  $\widetilde{\mathcal{N}} =  \cup_{j \in \mathcal{N}}\mathcal{N}_j$ represents the set of active units  within the power network and $\widetilde{\mathcal{E}} \subseteq \widetilde{\mathcal{N}} \times \widetilde{\mathcal{N}}$ the set of connections. 
Moreover, we let  $H  \in \mathbb{R}^{|\widetilde{\mathcal{N}}| \times |\widetilde{\mathcal{E}}|}$ 
be the incidence matrix of ($\widetilde{\mathcal{N}}, \widetilde{\mathcal{E}}$).
In addition, the following variables are defined for compactness in presentation,
\begin{subequations}\label{eq_s}
\begin{align}
    s_j^T &= [(-p^M_j)^T \;,\; (d^c_j)^T], j \in \mathcal{N}, \\
    \widetilde{s}_j &= s_j + p^L_j, j \in \mathcal{N},
\end{align}
\end{subequations}
where $\widetilde{s}  \in \mathbb{R}^{|\widetilde{\mathcal{N}}|}$  is a vector with all generation and controllable and uncontrollable demand units.

The proposed \textit{Extended Primal-Dual scheme}, is presented below
\begin{subequations}\label{eq5}
\begin{align}
    \widetilde{\Gamma} \dot{\psi} &= H^T p^c
    \label{5a},
    \\
    \Gamma \dot{p}^c &= \widetilde{s} - H \psi,
    \label{5b} \\
    u  &= p^c, \label{5c}
    \end{align}
\end{subequations}
where $\widetilde{\Gamma} \in \mathbb{R}^{|\widetilde{\mathcal{E}}| \times |{\widetilde{\mathcal{E}}|}}$ and  $\Gamma \in \mathbb{R}^{|\widetilde{\mathcal{N}}| \times |\widetilde{\mathcal{N}}|}$ are diagonal matrices containing the positive time constants associated with \eqref{5a} and \eqref{5b} respectively, and  $p^c_{k,j}$ corresponds to the power command variable associated with active unit   $k$ at bus $j$, that is also used as the input to \eqref{eq2} following \eqref{5c}.
A schematic representation of the system \eqref{eq1}, \eqref{eq2}, \eqref{eq_s},  \eqref{eq5} is provided in Fig. \ref{Fig1}.

{
\begin{Remark}\label{rem_comm_network}
The proposed \textit{Extended Primal-dual scheme} assumes communication among prosumption units by considering the connected graph ($\widetilde{\mathcal{N}}, \widetilde{\mathcal{E}}$). 
Note that when the communication of prosumption is considered for the Primal-Dual scheme, then its communication topology (i.e.  a meshed network at bus level with a star structure within each bus to {enable} communication from prosumption units towards the bus  controller) is a special case to that of the \textit{Extended Primal-dual scheme}.
It should also be noted that, apart from connectivity, no assumption is made on the topology of the communication network. 
The latter allows practical {aspects} to be considered in the design of the communication network (e.g. communication among buses could be at bus level only).
\end{Remark}
}

Following the \textit{Extended Primal-Dual scheme}, privacy-seeking users share power command signals instead of their generation and demand values.
Hence,  the prosumption profiles are not communicated towards local controllers.
The latter suffices to ensure  privacy against naive eavesdroppers, {since inferring the prosumption profiles from the power command variables would require knowledge of the underlying dynamics.
The privacy of prosumption profiles against naive eavesdroppers under the \textit{Extended Primal-Dual scheme} is demonstrated in the following proposition.

\begin{proposition}\label{prop_privacy_naive}
Consider any supply unit $k,j$ implementing the \textit{Extended Primal-Dual scheme} \eqref{eq5}. 
Then, its prosumption profile $\widetilde{s}_{k,j}$ is private against eavesdroppers with knowledge of K1. 
\end{proposition}
}

{
\begin{Remark}\label{rem_Extended}
Proposition \ref{prop_privacy_naive} demonstrates that the \textit{Extended Primal-Dual scheme}  ensures the privacy of prosumption profiles against naive eavesdroppers.
In addition, the \textit{Extended Primal-Dual scheme} achieves the same objectives as the \textit{Primal-Dual scheme}, both being inspired by suitable dual decomposition approaches. 
In particular,  it ensures generation/demand balance and that the frequency attains its nominal value at steady state.
The latter follows by summing \eqref{eq1b} and \eqref{5b} at steady state over all units in $\mathcal{N}$ and $\widetilde{\mathcal{N}}$, which yield
$\sum_{j \in \mathcal{N}} \zeta^*_j - D_j \omega^*_j = 0$ and
 $\sum_{j \in \widetilde{\mathcal{N}}} \widetilde{s}^* = 0$ respectively, allowing to deduce that $\omega^* = 0$ due to \eqref{eq1a}.
 The second objective is to ensure that power command variables are synchronized at steady state, which is achieved through \eqref{5a}.
 Hence, the \textit{Extended Primal-Dual scheme} achieves the objectives of the \textit{Primal-Dual scheme} and  additionally ensures privacy against naive eavesdroppers.
\end{Remark}
}

\subsection{Equilibrium Analysis}

We now provide a definition of an equilibrium point to the interconnected dynamical system \eqref{eq1}, \eqref{eq2}, \eqref{eq_s},  \eqref{eq5}. 

\begin{definition} 
The point $\alpha^*$ = $(\eta^*,\psi^*,\omega^*,x^*, p^{c,*})$ defines an equilibrium of the system \eqref{eq1}, \eqref{eq2}, \eqref{eq_s},  \eqref{eq5} if all time derivatives of \eqref{eq1}, \eqref{eq2}, \eqref{eq_s},  \eqref{eq5} are equal to zero at this point.
\end{definition} 

We will make use of the following equilibrium equations for \eqref{eq1}, \eqref{eq2}, \eqref{eq_s},  \eqref{eq5}.
\begin{subequations}\label{eq7}
\begin{align}
    &0=\omega^*_{i}-\omega^*_{j}, (i,j)\in{\mathcal{E}}, \label{eq7a} \\
  \begin{split}
&0  = \vect{1}^T_{|N^G_j|} p^{M,*}_{j} - \vect{1}^T_{|N^L_j|}  d^{c,*}_{j} -\vect{1}^T_{|N_j|} p^L_{j} \\
 & \hspace{2mm} -  \sum_{i \in \mathcal{N}^s_j} p^*_{ji} +
  \sum_{i \in \mathcal{N}^p_j} p^*_{ij}, j \in \mathcal{N} \label{eq7b}
\end{split} \\
&0 = -x^*_{k,j} + m_{k,j}(u^*_{k,j} - \omega^*_j), k \in \mathcal{N}^G_j, j\in{\mathcal{N}}, \label{eq7c}
\\
 &0 = H^T p^{c,*},
  \label{7f}\\
   &0 = \widetilde{s}^* - H \psi^*,
    \label{7g}
\\
\text{where} & \text{ the variables $p^*, p^{M,*}, d^{c,*}, u^*, \widetilde{s}^*$ satisfy} \nonumber\\
p^*_{ij} &=B_{ij}\eta^*_{ij}, (i,j)\in{\mathcal{E}}, \label{eq7h} \\
 p^{M,*}_{k,j} &= x^*_{k,j} + h_{k,j} (u^*_{k,j} - \omega^*_j), k \in \mathcal{N}^G_j, j\in{\mathcal{N}}, \label{eq7d}
 \\
 d^{c,*}_{k,j} &= - h_{k,j} (u^*_{k,j} - \omega^*_j), k \in \mathcal{N}^L_j, j \in{\mathcal{N}}, \label{eq7e} \\
 u^*_{k,j} &= p^{c,*}_{k,j},  k \in \mathcal{N}_j, j \in{\mathcal{N}} \label{eq7i}, \\
(s_j^*)^T &= [(-p^{M,*}_j)^T \;,\; (d^{c,*}_j)^T], j \in \mathcal{N}, \\
\widetilde{s}^* &= s^* + p^L.
 \end{align}
\end{subequations}
The following lemma characterizes the equilibria of \eqref{eq1}, \eqref{eq2}, \eqref{eq_s},  \eqref{eq5}.

\begin{lemma}\label{Lemma1}
The equilibria of \eqref{eq1}, \eqref{eq2}, \eqref{eq_s},  \eqref{eq5} satisfy $\omega^*=\vect{0}_{|{\mathcal{N}}|}$ and
$p^{c,*} \in \Ima(\vect{1}_{|\widetilde{\mathcal{N}}|})$.
\end{lemma}

Lemma \ref{Lemma1} demonstrates that the presented scheme ensures that the frequency attains its nominal value at equilibrium, which is a main objective of secondary frequency control. In addition, it shows that power command variables share the same value at steady state. The latter can be used to enable an optimal power allocation, as demonstrated in the following section.

\subsection{Optimality and Stability Analysis}

The following proposition  provides necessary and sufficient conditions that ensure that the equilibrium values of  $p^{M}$ and $d^{c}$ are global solutions to the optimization problem \eqref{eq3}.

\begin{proposition}\label{proposition1}
Let $q_{k,j}(m_{k,j}+h_{k,j})=1, k \in \mathcal{N}^G_j, j \in \mathcal{N}$ and $q_{k,j}h_{k,j}=1, k\in \mathcal{N}^L_j, j \in \mathcal{N}$. 
Then, the equilibrium values $p^{M,*}$ and $d^{c,*}$ of system \eqref{eq1}, \eqref{eq2}, \eqref{eq_s},  \eqref{eq5} globally minimize {the} optimization problem \eqref{eq3}.
\end{proposition}

Proposition \ref{proposition1} follows directly from the KKT conditions \cite{boyd2004convex}. It demonstrates how the controller gains in generation and controllable load units should be designed such that an optimal power allocation is ensured.
Hence, we deduce that the \textit{Extended Primal-Dual scheme} \eqref{eq5} enables an optimal power allocation.

The following theorem provides global asymptotic stability guarantees for \eqref{eq1}, \eqref{eq2}, \eqref{eq_s},  \eqref{eq5}.

\begin{theorem}\label{theorem1}
Solutions to \eqref{eq1}, \eqref{eq2}, \eqref{eq_s},  \eqref{eq5} globally asymptotically converge to the set of its equilibria, where $\omega^*= \boldsymbol{0}_{\mathcal{|N|}}$.
\end{theorem}

Theorem \ref{theorem1} guarantees the convergence of solutions to \eqref{eq1}, \eqref{eq2}, \eqref{eq_s},  \eqref{eq5} to the set of its equilibria. 
In addition,  the \textit{Extended Primal-Dual scheme} is {locally} verifiable and applicable to arbitrary network configurations.
Furthermore,  the presented scheme guarantees the privacy  of the prosumption profiles against naive eavesdroppers.
Noting also that Proposition \ref{proposition1} demonstrates how optimality may be achieved at steady state, it follows that the \textit{Extended Primal-Dual scheme} satisfies all objectives of Problem \ref{problem_definition}, except from ensuring privacy against intelligent eavesdroppers.

\subsection{Discussion}

The scheme presented in this section extends the \textit{Primal-Dual scheme} \eqref{eq4} by including a controller at each unit contributing to secondary frequency control.
The \textit{Extended Primal-Dual scheme} results in the transmission of power command signals instead of prosumption signals, 
which enables privacy against naive eavesdroppers, {as demonstrated by Proposition \ref{prop_privacy_naive}}.
{On the other hand, the interaction between an increased number of controllers may result in slower convergence.}
The proposed scheme yields an optimal power allocation, ensures that frequency attains its nominal value at steady state and guarantees the global stability of the power network as follow from Proposition \ref{proposition1},  Lemma \ref{Lemma1} and Theorem \ref{theorem1} respectively.
However, the \textit{Extended Primal-Dual scheme} \eqref{eq5} does not ensure the privacy of generation and demand profiles against intelligent eavesdroppers.
In particular, an intelligent eavesdropper may use the communicated power command trajectories and knowledge of the underlying power command dynamics to infer the  prosumption profiles
{by reversing \eqref{5b}, i.e. using} $\widetilde{s} =  \Gamma \dot{p}^c + H \psi$. 
In the next section, we present {a} scheme that aims to resolve this issue.

\section{Privacy-Preserving scheme}\label{IV}

In this section we present {a} scheme
that aims to  preserve the beneficial properties of the \textit{Extended Primal-Dual scheme} described in the previous section and  simultaneously {guarantee the} privacy of the generation/demand profiles against intelligent eavesdroppers.

\subsection{Privacy-Preserving scheme}
The proposed scheme, which shall be referred {to} as the \textit{Privacy-Preserving scheme}, incorporates a privacy-enhancing signal $n$ in the power command dynamics, as follows
\begin{subequations}\label{eq8}
\begin{align}
    \widetilde{\Gamma} \dot{\psi} &= H^T p^c
    \label{eq8a},
    \\
\Gamma \dot{p}^c &= \widetilde{s} - H \psi + n
  \label{eq8b}, \\
  u &= p^c. \label{eq8c}
\end{align}
\end{subequations}
In \eqref{eq8} above, the locally Lipschitz, privacy-enhancing signal $n = [n_i]_{i \in {\mathcal{N}}}$, where $n_i = [n_{k,i}]_{{k \in {\mathcal{N}_i}}}$,  adapts the derivative of the power command variables 
to enable enhanced privacy properties.

The design of the signal $n$ is crucial in providing enhanced privacy properties and simultaneously enabling stability and optimality guarantees for the \textit{Privacy-Preserving scheme} \eqref{eq8}.
Some desired properties of the privacy-enhancing signal $n$ are: (i) to permit the existence of equilibria, by taking a constant value when the states of the system are at equilibrium, and
(ii)~to enable an optimality interpretation of the resulting equilibria.
Both objectives can be achieved if $n$ is zero at steady state since in this case the equilibria of \eqref{eq1}, \eqref{eq2},  \eqref{eq_s}, \eqref{eq8},  and \eqref{eq1}, \eqref{eq2}, \eqref{eq_s},  \eqref{eq5} are identical.

The following design condition is imposed on the privacy-enhancing signal $n$.
As demonstrated below, this condition ensures the privacy of the prosumption profiles against intelligent eavesdroppers  and 
allows stability and optimality to be deduced.
{It should be noted that the trajectories of $n$ are in general non-unique.}

\begin{design}\label{assumption1}
{The privacy-enhancing signals  satisfy $n_{k,j} = n^d_{k,j} + n^f_{k,j}, k \in \mathcal{N}_j, j \in \mathcal{N}$, where:
\begin{enumerate}[(i)]
    \item $n^d_{k,j}(t) = - \xi_{k,j}(t) \dot{p}^c_{k,j}(t)$,  where 
     the non-negative signal $\xi_{k,j}(t)$ satisfies $\dot{\xi}_{k,j}(t) < \hat{\beta}_{k,j}$ for all $t \geq 0$,
    \item $|n^f_{k,j}(t)| < \beta_{k,j}|\omega_j(t)|$, 
    for all $t \geq 0$.
\end{enumerate}
   Moreover,   the positive design constants $\beta_{k,j}, \hat{\beta}_{k,j}$ satisfy
   $\begin{bmatrix} -h_{k,j} -D_j/|\mathcal{N}_j| & h_{k,j} + \beta_{k,j}/2\\
     h_{k,j} + \beta_{k,j}/2 & -h_{k,j} + \hat{\beta}_{k,j}/2
    \end{bmatrix} \preceq 0, k \in \mathcal{N}_j, j \in \mathcal{N}$.}
\end{design}

Design Condition \ref{assumption1} splits the privacy-enhancing signal $n$ to two other signals,  $n^d$ and $n^f$, that serve different purposes.
The signal $n^d_{k,j}$ is proportional to the power command derivative $\dot{p}^c_{k,j}$ with {a non-negative, time-varying}  gain $\xi_{k,j}$   {designed such that  $\dot{\xi}_{k,j}(t) < \hat{\beta}_{k,j}$ is satisfied at all times}. 
The latter adjusts the rate at which the power command variables respond to external signals and  makes any prior estimates of the power command model inaccurate.
Hence, a potential eavesdropper utilizing model-based observations will produce inaccurate results. 
The component $n^f$ introduces a noise signal\footnote{It should be noted that    $n^f$ {(and similarly $\xi$)} are treated as time-dependent variables rather than random variables,
following the assumption that $n$ is locally Lipschitz.
The latter is made for simplicity and to avoid a diversion of the paper focus from the privacy properties of the proposed schemes.} 
 that is mixed with the generation/demand values.
The latter offers improved privacy properties since: (i) the generation/demand profile information in the controller is distorted, and (ii) it perturbs the communicated signals of all controllers when a disturbance occurs,  making it harder to detect the origin of the disturbance from a change in the transmitted signal.
Design Condition \ref{assumption1}(ii)
restricts the magnitude of $n^f$ in relation with the magnitude of the local frequency. 
The {values of $\beta_{k,j}, \hat{\beta}_{k,j}$ are selected to satisfy the linear matrix inequality (LMI)  in Design Condition \ref{assumption1}} such that convergence is guaranteed, as demonstrated in Theorem \ref{theorem2} later on. 
  These properties enable the privacy of prosumpion against intelligent eavesdroppers since the same power command trajectories result from a (wide) class of prosumption profiles due to different potential trajectories of the privacy-enhancing signal $n$.
  {The latter is analytically demonstrated in Section \ref{sec_privacy_analysis} below.}
In addition, note that  since all communicated power command signals synchronize at steady state, 
{their equilibrium values do not convey any information about local generation/demand.}

{
\begin{Remark}
The bounds $\beta_{k,j}$ and $\hat{\beta}_{k,j}$ associated with $n^f_{k,j}$ and $n^d_{k,j}$ respectively are interdependent through the LMI  in Design Condition \ref{assumption1}.
Hence, there is a trade-off between the maximum allowed derivative of the gain $\xi_{k,j}$ and the maximum magnitude ratio  between the signal $n^f_{k,j}$ and the local frequency $\omega_j$.
The latter can be used for design purposes by placing different weights on the the associated bounds, and hence the effect, of signals  $n^f_{k,j}$ and $n^d_{k,j}$.
\end{Remark}
}

{
\begin{Remark}\label{rem_controllers}
The implementation of the proposed \textit{Extended Primal-Dual} and \textit{Privacy-Preserving schemes} requires suitable monitoring and communication capabilities and an increased number of controllers. 
These requirements are facilitated by  the significant improvement in monitoring, control and communication technologies  and their implementation on smart power grids \cite{pliatsios2020survey,   yang2011communication}.
It should also be noted that privacy issues are mostly motivated as side effects of  these improvements  which in most cases coincide with enhanced control capabilities.
\end{Remark}
}

{
\subsection{Privacy analysis}\label{sec_privacy_analysis}

In this section, we present our main privacy results regarding the proposed \textit{Privacy-Preserving scheme}.
First, we clarify that for an intelligent eavesdropper, K1 implies knowledge of all power command signals communicated to and from a considered unit. In addition, K2 implies  knowledge of the 
\emph{Privacy-Preserving scheme} dynamics \eqref{eq8}.

The following proposition demonstrates that the proposed scheme preserves the privacy of the prosumption profiles against intelligent eavesdroppers.

\begin{proposition}\label{prop_privacy}
Consider any supply unit $k,j$ implementing the \textit{Privacy-Preserving scheme} \eqref{eq8}. 
Then, its prosumption profile $\widetilde{s}_{k,j}$ is private against intelligent eavesdroppers with knowledge of K1 and K2.
\end{proposition}

Proposition \ref{prop_privacy} provides privacy guarantees for the prosumption profiles when the \textit{Privacy-Preserving} scheme is implemented.
The latter demonstrates that the proposed scheme satisfies objective (i) within Problem \ref{problem_definition}.

A reasonable case to be considered is when intelligent eavesdroppers gain knowledge of the steady-state value of the privacy-preserving signal $n$, i.e. have the following knowledge:   \newline
(K3) The steady state value of the privacy-preserving signal $n$, i.e. that $\lim_{t \rightarrow \infty} n(t) = \vect{0}_{|\widetilde{\mathcal{N}}|}$.

The following proposition shows that the prosumption trajectories are private against eavesdroppers with knowledge of K1, K2 and K3.
We remind that the definition of a private trajectory is provided in Definition \ref{dfn_privacy}(i).

\begin{proposition}\label{prop_privacy_trajectory}
Consider any supply unit $k,j$ implementing the \textit{Privacy-Preserving scheme} \eqref{eq8}. 
Then, its prosumption trajectory is private against intelligent eavesdroppers with knowledge of K1, K2 and K3.
\end{proposition}

Note that Proposition \ref{prop_privacy_trajectory} does not guarantee the privacy of the prosumption at steady state, since knowledge of the variable $\psi$ may yield the equilibrium values of $\widetilde{s}$ from \eqref{7g}.
However, since $\psi$ results from integrating the differences between communicated power command variables, any inaccuracy on determining these variables will lead to growing deviations between the estimated and true values of $\psi$, compromising the reliability of such estimate.

Stronger privacy guarantees may be obtained, such that the prosumption profiles are kept private when eavesdroppers have knowledge of K3, by relaxing K1.
In particular, we consider the case where an eavesdropper does not have full knowledge of the information communicated to a considered unit, i.e. has knowledge of:  \newline
(K4) A proper subset of the power command signals  communicated to and from a given unit, for which it aims to obtain private information. 

The following proposition guarantees the privacy of prosumption profiles when intelligent eavesdroppers have knowledge of K2, K3 and K4.

\begin{proposition}\label{lemma_privacy}
Consider any supply unit $k,j$ implementing the \textit{Privacy-Preserving scheme} \eqref{eq8}. 
Then, its prosumption profile $\widetilde{s}_{k,j}$ is private against intelligent eavesdroppers with knowledge of K2, K3 and K4.
\end{proposition}

Proposition \ref{lemma_privacy} enables privacy guarantees of the prosumption profile, when knowledge of the steady state value of $n$ is available. However, it assumes  that the intelligent eavesdropper does not possess full knowledge of the information communicated to and from the considered privacy-seeking unit. Note that the latter case might describe eavesdroppers associated with some  prosumption unit that communicates with the considered privacy-seeking unit, under specific conditions on the communication network topology such that K4 is satisfied.

\begin{Remark}\label{rem_noise_signals}
The presented privacy results hold when either of $n^f$ or $n^d$ is neglected in $n$.
However, their combined impact keeps additional information associated with the prosumption profiles private.
In particular, the presence of $n^f$ results in a change on all controllers after a power disturbance, making it difficult to infer its origin from the power command signals.
In addition, $n^d$ makes model based inference inaccurate, and hence difficult to have a reasonable range estimate of the prosumption magnitude, i.e. obtain an estimate with  a margin of error analogous to the magnitude of $n^f_{k,j}$.
\end{Remark}
}

\subsection{Optimality and Stability Analysis}

In this section we provide analytic optimality and stability guarantees for system \eqref{eq1}, \eqref{eq2},  \eqref{eq_s}, \eqref{eq8}.

The following proposition  extends Proposition \ref{proposition1} by demonstrating that Design Condition \ref{assumption1} enables an optimal steady state power allocation.

\begin{proposition}\label{proposition2}
Let Design Condition \ref{assumption1},  $q_{k,j}(m_{k,j}+h_{k,j})=1, k \in \mathcal{N}^G_j, j \in \mathcal{N}$ and $q_{k,j}h_{k,j}=1, k\in \mathcal{N}^L_j, j \in \mathcal{N}$ hold. Then, the equilibrium values  $p^{M,*}$ and $d^{c,*}$ of system \eqref{eq1}, \eqref{eq2},  \eqref{eq_s}, \eqref{eq8},  globally minimize the optimization problem \eqref{eq3}.
\end{proposition}

Proposition \ref{proposition2} demonstrates that when Design Condition \ref{assumption1}, and the  gain conditions provided in Proposition \ref{proposition1} hold, then the \textit{Privacy-Preserving scheme} yields an optimal power allocation.
The latter follows trivially from Proposition \ref{proposition1}, since  the privacy-enhancing signal $n_{k,j}$ is zero at steady state from Design Condition \ref{assumption1}, which results in identical equilibrium points for \eqref{eq1}, \eqref{eq2},  \eqref{eq_s}, \eqref{eq8} and \eqref{eq1}, \eqref{eq2}, \eqref{eq_s},  \eqref{eq5}.

The following theorem  demonstrates that when Design Condition \ref{assumption1} holds, then the set of equilibria of \eqref{eq1}, \eqref{eq2},  \eqref{eq_s}, \eqref{eq8},  is attracting. The latter  shows that the proposed  \textit{Privacy-Preserving scheme} does not compromise the stability of the power network.

\begin{theorem}\label{theorem2}
Let Design Condition \ref{assumption1} hold. Then, the solutions of \eqref{eq1}, \eqref{eq2},  \eqref{eq_s}, \eqref{eq8},  globally asymptotically converge to the set of its equilibria,  where $\omega^*= \vect{0}_{\mathcal{|N|}}$.
\end{theorem}

Theorem \ref{theorem2} guarantees the convergence of solutions to \eqref{eq1}, \eqref{eq2},  \eqref{eq_s}, \eqref{eq8},  to the set of its equilibria.
In addition, the dynamics of \eqref{eq1}, \eqref{eq2},  \eqref{eq_s}, \eqref{eq8},  are distributed, applicable to arbitrary network configurations and locally verifiable.
    Moreover, {as demonstrated in Section \ref{sec_privacy_analysis},  the}  \textit{Privacy-Preserving scheme} enables the privacy of prosumption profiles against informed eavesdroppers.
Finally, as demonstrated in Proposition \ref{proposition2}, the presented scheme enables an optimal power allocation among generation and controllable demand.
Hence, all objectives of Problem \ref{problem_definition} are satisfied.

\section{Simulation on the NPCC 140-bus system}\label{V}

In this section, we {validate} our analytic results with simulations using the Power system toolbox \cite{chow1992toolbox} on Matlab. 

\subsection{{Test system description}}\label{sec_test_system}

We use the Northeast Power Coordinating Council (NPCC) 140-bus interconnection system. This model is more detailed and realistic than the considered analytical model, including voltage dynamics,  line resistances, and a transient reactance generator model\footnote{{The details of the simulation model can be found in the Power System Toolbox data file datanp48.}}.
The test system consists of 93 load buses 
and 47 generation buses and has a total real power of 28.55 GW.

\subsection{{Simulation Results}}

{For our simulations,} controllable demand was considered in $20$ load buses, where at each bus the number of controllable loads was randomly selected from an integer uniform distribution with range $[90, 180]$.
A single generation unit was added at each of $20$ generation buses.
In addition, quadratic cost functions were considered for generation and controllable demand following the description in \eqref{eq3}.
The values for $q_{k,j}, k \in \mathcal{N}_j, j \in \mathcal{N}$ were selected from a uniform distribution with range $[50, 250]$.
For the simulation, a step change in demand of magnitude $0.2$ per unit (100 MW) at $10$ randomly selected loads at each of buses $2$ and $3$ was considered at $t = 1$ second.
The time step for the simulations, {denoted by $\Delta T$,} was set at 10 ms.

 \begin{figure}[t!]
\centering
\includegraphics[scale = 0.57]{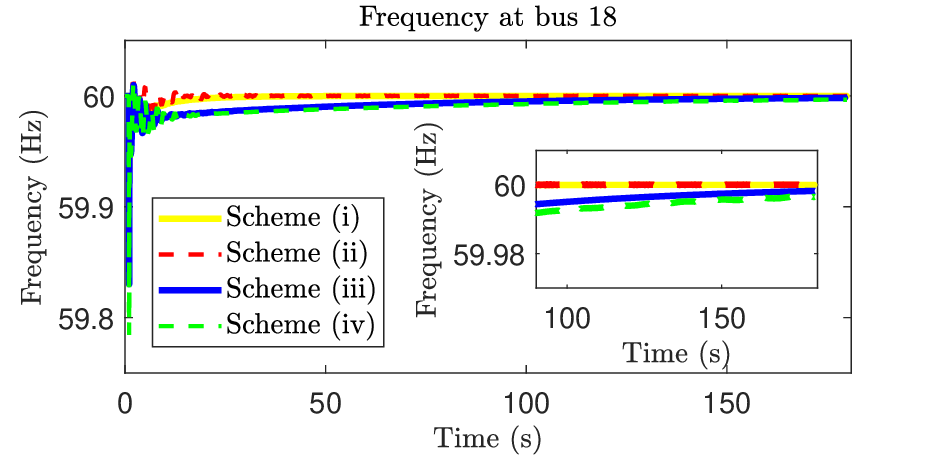}
\vspace{-5mm}
\caption{Frequency at bus $18$ when the following  schemes are implemented: (i) Integral action scheme, (ii) Primal-Dual scheme, (iii) Extended Primal-Dual scheme, and (iv) Privacy-Preserving scheme.}
\label{Fig_frequency}
\vspace{-5mm}
\end{figure}

\begin{figure*}[!htb]
\begin{minipage}{0.9\textwidth}
\begin{subfigure}{0.23\textwidth}
\centering
\includegraphics[scale=0.5]{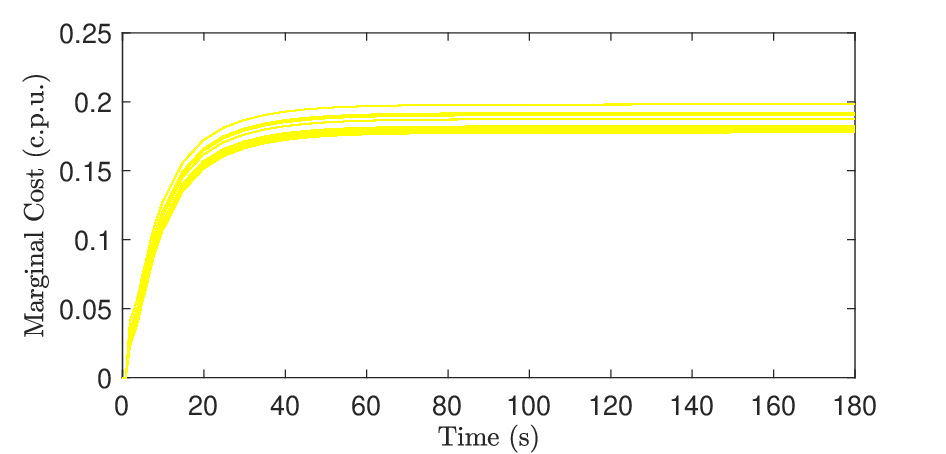}
\end{subfigure}
\hspace{0.24\textwidth}
\begin{subfigure}{0.24\textwidth}
\centering
\includegraphics[scale=0.5]{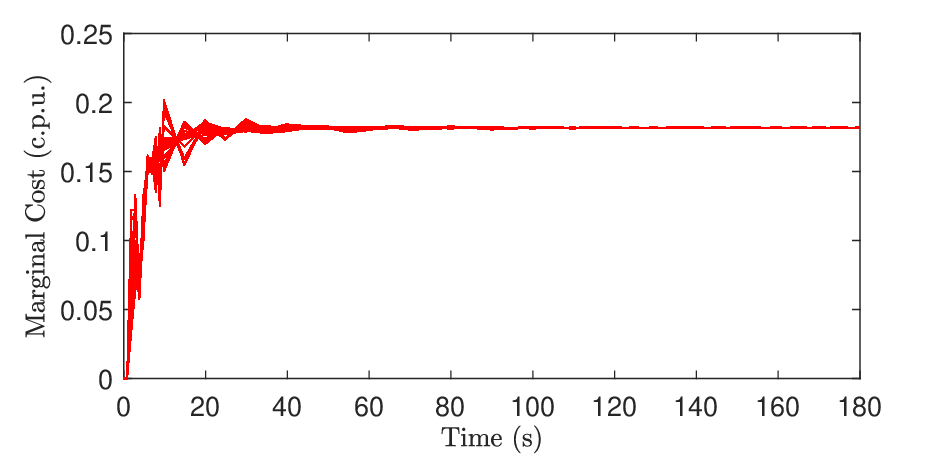}
\end{subfigure}
\\
\begin{subfigure}{0.23\textwidth}
\includegraphics[trim = 0 0 0  2mm, clip,scale=0.5]{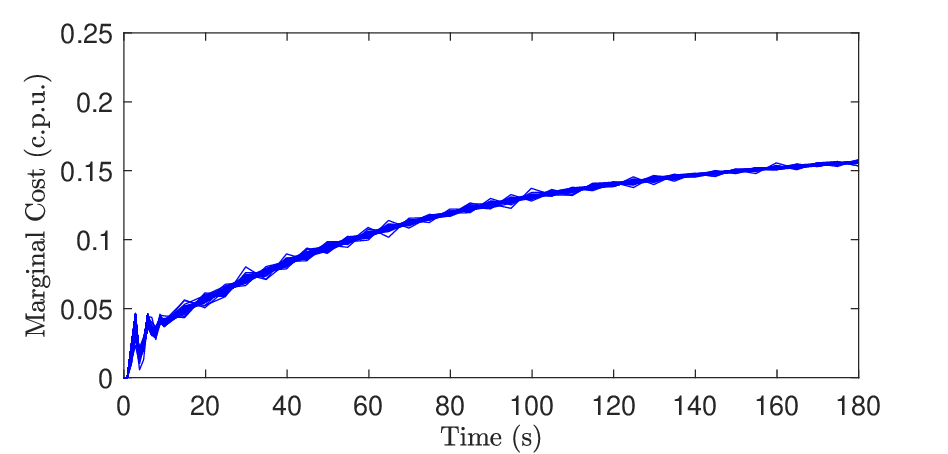}
  \end{subfigure} 
\hspace{0.24\textwidth}
\begin{subfigure}{0.24\textwidth}
\centering
 \includegraphics[trim = 0 0 0  2mm, clip,scale=0.50]{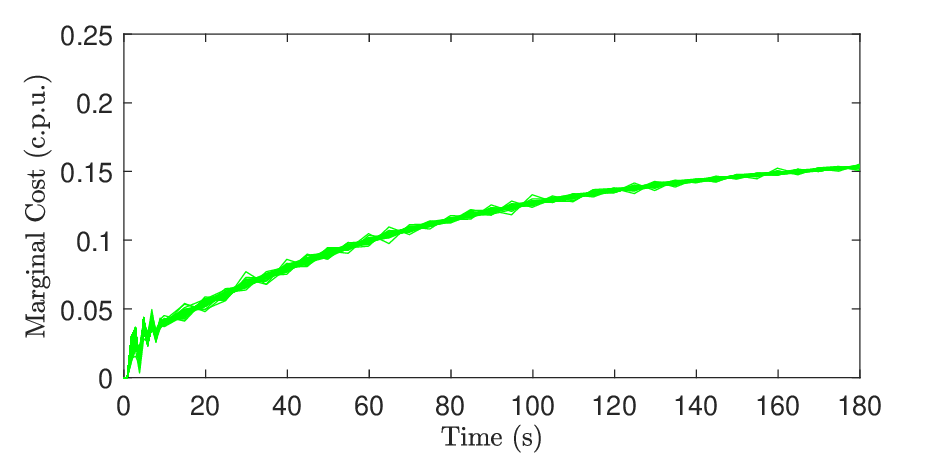}
\end{subfigure}
\end{minipage}
 \hspace{-7.5mm} 
\begin{minipage}{0.08\textwidth}
{\includegraphics[trim = 77mm 15mm 17mm 12mm, clip,scale=0.62]{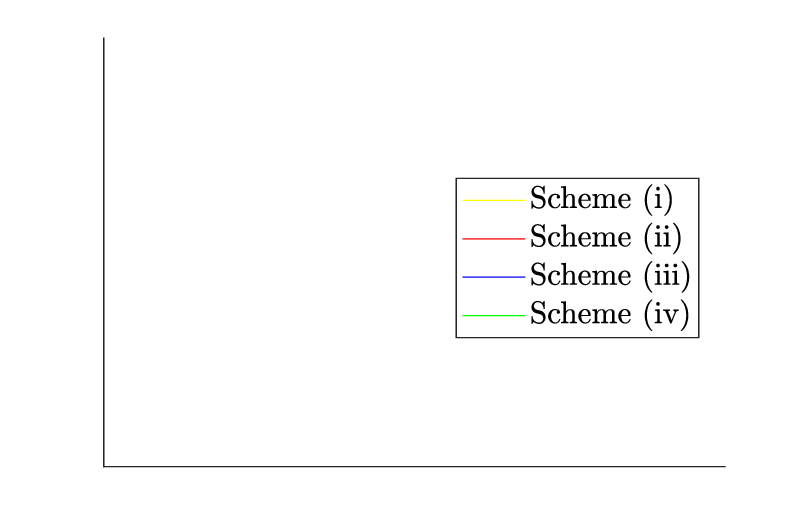}}
 \end{minipage}
 \vspace{-1mm}
\caption{Marginal costs for all generation and controllable demand units contributing to secondary frequency control when the following  control schemes are implemented: (i) Integral action scheme, (ii) Primal-Dual scheme, (iii) Extended Primal-Dual scheme, and (iv) Privacy-Preserving scheme.}
\label{Fig_marginal}
\vspace{-2mm}
\end{figure*}

The system was tested {under} the four control schemes described below:
\begin{enumerate}[(i)]
    \item An Integral action scheme, where generation units and controllable loads integrate the local frequency with the controller gains selected to be inversely proportional to their respective cost coefficients.
    \item The \textit{Primal-Dual scheme}, described by \eqref{eq4}.
    \item  The \textit{Extended Primal-Dual scheme} that we proposed, described by \eqref{eq5}.
    \item  The \textit{Privacy-Preserving scheme} that we proposed, described by \eqref{eq8} and Design Condition \ref{assumption1}. 
    {First suitable values for $\beta_{k,j}, \hat{\beta}_{k,j}, k \in \mathcal{N}_j, j \in \mathcal{N}$ were selected in accordance with the LMI in Design Condition \ref{assumption1}. 
    The values of $\xi_{k,j}(t)$  were then randomly selected at each time step such that 
    $(\xi_{k,j}(t) - \xi_{k,j}(t- \Delta T))/\Delta T$
    lied in   $[- \hat{\beta}_{k,j},  \hat{\beta}_{k,j}]$ following  Design Condition \ref{assumption1}(i).
    In addition, the values of $n^f_{k,j}(t)$ were randomly selected at each time step from the uniform distribution $[-\beta_{k,j}|\omega_j|,\beta_{k,j}|\omega_j|]$ such that  Design Condition \ref{assumption1}(ii) was satisfied}. 
\end{enumerate}
In schemes (ii)-(iv), the dynamics of the implemented generation and controllable demand units followed from \eqref{eq2} and the controller gains were selected such that the optimality conditions presented in Propositions \ref{proposition1} and \ref{proposition2} were satisfied.
The communication network associated with scheme (ii) had the same structure as the power network.
A random connected communication network was generated when schemes (iii) and (iv) were implemented.
For consistency, the same sets of randomly selected parameters were considered in all simulations.

\begin{figure*}[ht]
\begin{subfigure}{0.23\textwidth}
\centering
\includegraphics[scale=0.4]{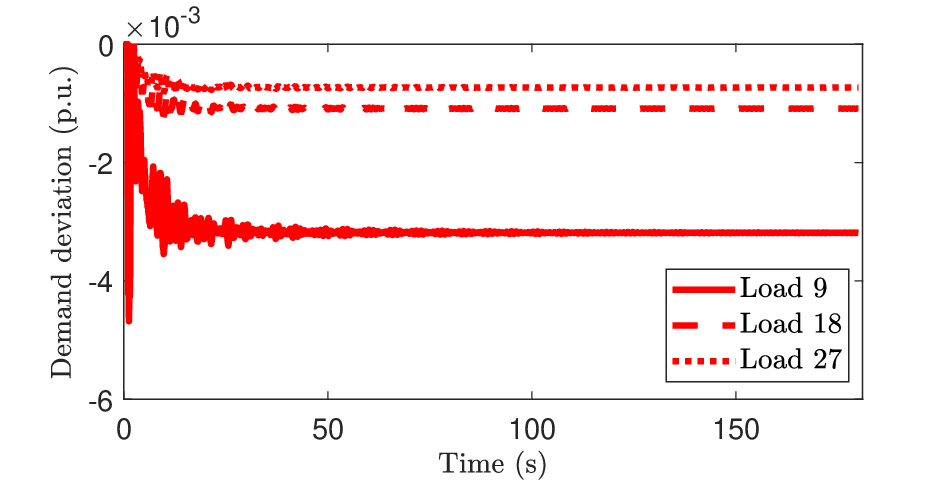}
\end{subfigure}
\hspace{0.09\textwidth}
\begin{subfigure}{0.24\textwidth}
\centering
\includegraphics[scale=0.4]{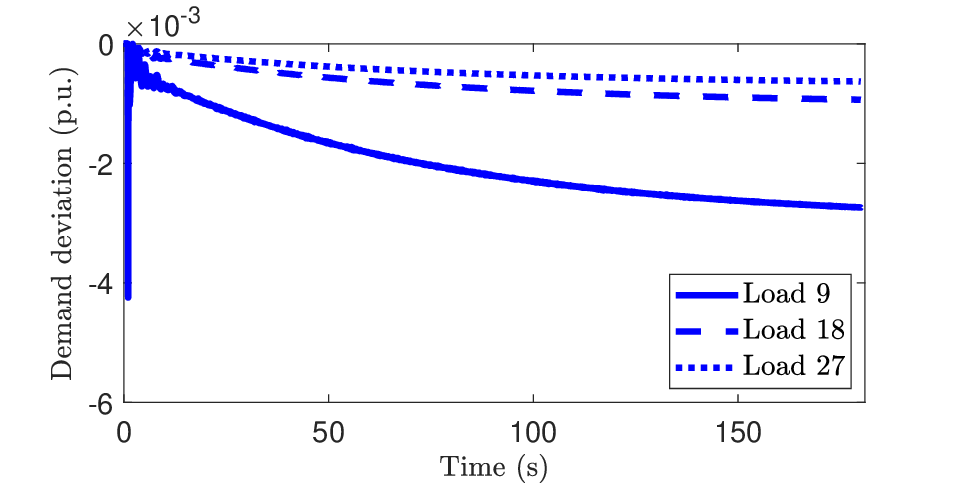}
\end{subfigure}
\hspace{0.09\textwidth}
\begin{subfigure}{0.24\textwidth}
\centering
\includegraphics[scale=0.4]{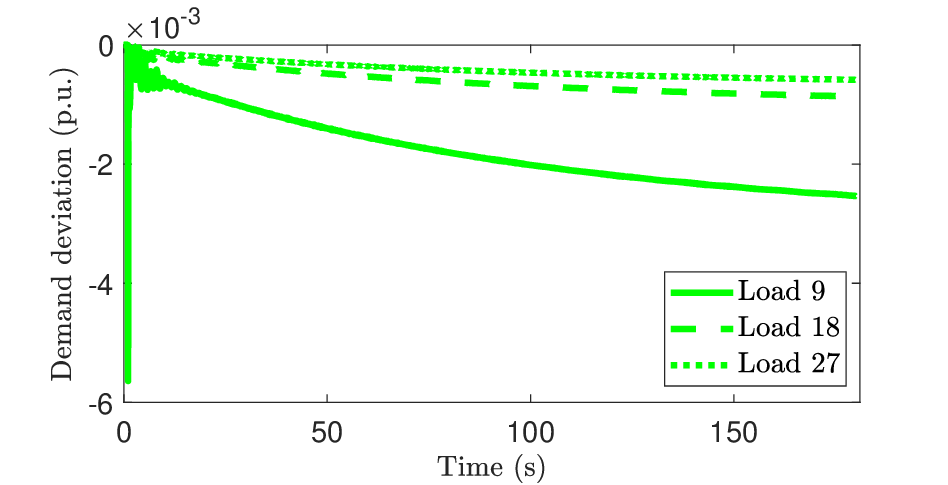}
\end{subfigure}
\\[1mm]
\begin{subfigure}{0.23\textwidth}
\includegraphics[trim = 0 0 0  2mm, clip,scale=0.4]{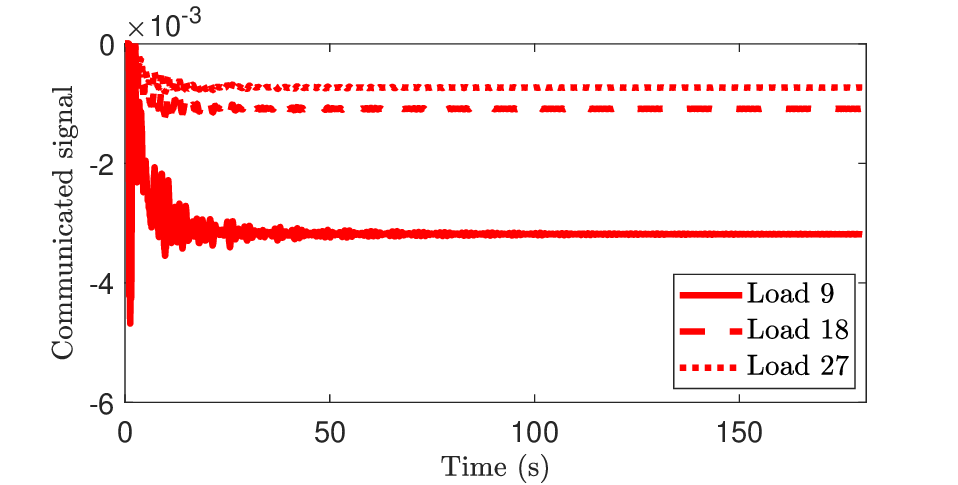}
  \end{subfigure} 
\hspace{0.09\textwidth}
\begin{subfigure}{0.24\textwidth}
\centering
 \includegraphics[trim = 0 0 0  2mm, clip,scale=0.4]{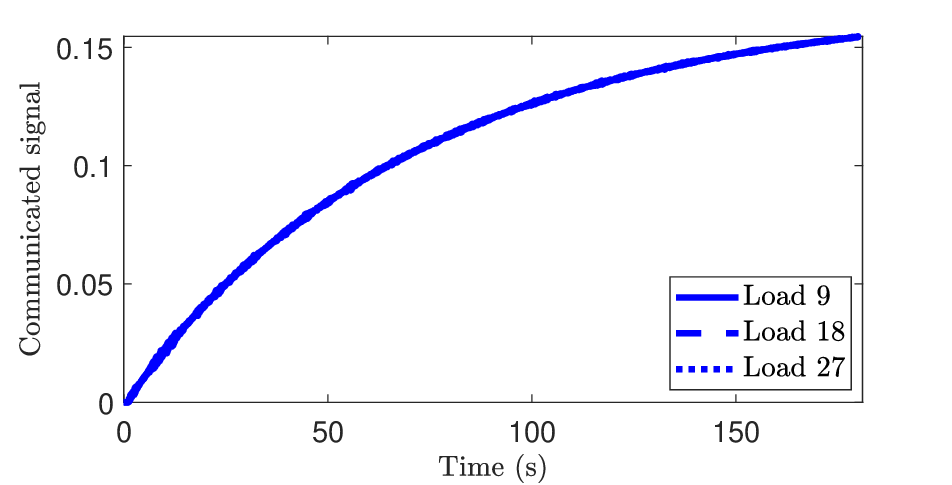}
\end{subfigure}
\hspace{0.09\textwidth}
\begin{subfigure}{0.24\textwidth}
\centering
 \includegraphics[trim = 0 0 0  2mm, clip,scale=0.4]{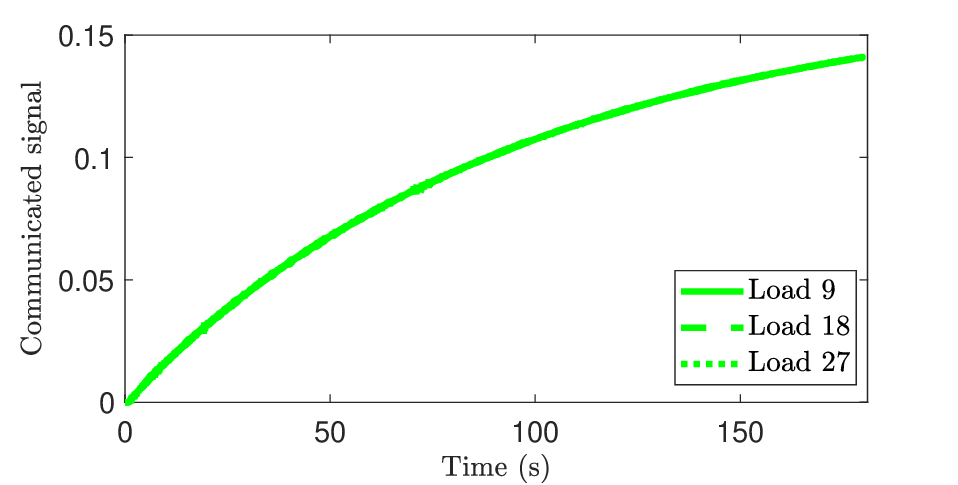}
\end{subfigure}
\vspace{-1mm}
\caption{Controllable demand (top) and communicated signals (bottom) for  loads $9, 18$ and $27$ at bus $2$  for the following control schemes: (left) Primal-Dual scheme, (middle) Extended Primal-Dual scheme, and (right) Privacy-Preserving scheme.}
\label{Privacy_A}
\vspace{-2mm}
\end{figure*}

The frequency response at a randomly selected bus (bus $18$) is depicted in Fig. \ref{Fig_frequency}.
From Fig. \ref{Fig_frequency}, it follows that the frequency converges to its nominal value at all simulated cases. The latter suggests that the proposed \textit{Extended Primal-Dual} and \textit{Privacy-Preserving} schemes yield a stable response. Note also that the frequency returns to within $0.01$ Hz from its nominal value in less than two minutes, which is well within the secondary frequency control timeframe.
Nevertheless, the \textit{Extended Primal-Dual}  and \textit{Privacy-Preserving} schemes result in slower convergence of frequency to its nominal value. {This is due to} a larger number of controllers {that} need to synchronize for convergence. In addition, the implementation of the \textit{Privacy-Preserving scheme}, and particularly  Design Condition \ref{assumption1}(i), 
results in slower convergence
compared with the \textit{Extended Primal-Dual scheme}.

To demonstrate the optimality of the proposed analysis, we consider the marginal costs of each active unit, defined as the absolute value of the cost derivative of the local cost functions.
The marginal costs  for all controllable loads and local generators{, measured in cost per additional unit (c.p.u.),} are depicted in Fig. \ref{Fig_marginal}. 
From Fig. \ref{Fig_marginal}, it follows that the marginal costs for all units converge to the same value {when  schemes (ii), (iii) and (iv) are employed}. The latter suggests that an optimal power allocation is attained at steady state and validates the presented optimality analysis. 
By contrast, the marginal costs differ at equilibrium in scheme (i), which suggests that a suboptimal response is obtained. 
{The latter follows since scheme (i) implements integral control action on the 
local frequency values.
Since local frequency trajectories are not identical, the integral action scheme will have a different impact at different buses, and hence yield a sub-optimal power sharing.}

To validate the enhanced privacy properties associated with the \textit{Extended Primal-Dual} and the \textit{Privacy-Preserving} schemes, compared with the \textit{Primal-Dual scheme}, we considered the communicated signals  from three randomly selected loads (loads $9, 18$ and $27$ at bus $2$).
The results are shown on Fig. \ref{Privacy_A}.
{Figure \ref{Privacy_A} demonstrates that the implementation of the \textit{Primal-Dual scheme} (scheme (ii)) compromises the privacy of the controllable demand units.
The latter follows since the demand values are communicated, as follows from the bottom-left subfigure in Fig. \ref{Privacy_A}, allowing a naive eavesdropper that possesses knowledge of the communicated information to infer the prosumption profile.}
By contrast, when the \textit{Extended Primal-Dual} and \textit{Privacy-Preserving} schemes were implemented (schemes (iii) and (iv) respectively), the privacy of the controllable load profiles against naive eavesdroppers is preserved {since, instead of prosumption profiles,  power command trajectories were communicated.}

\begin{figure}
\includegraphics[scale=0.57]{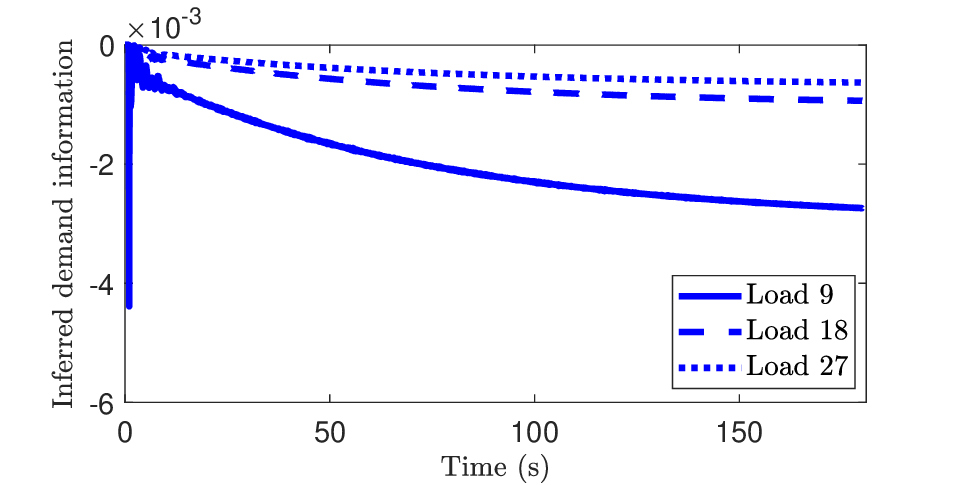}
\\
\includegraphics[scale=0.57]{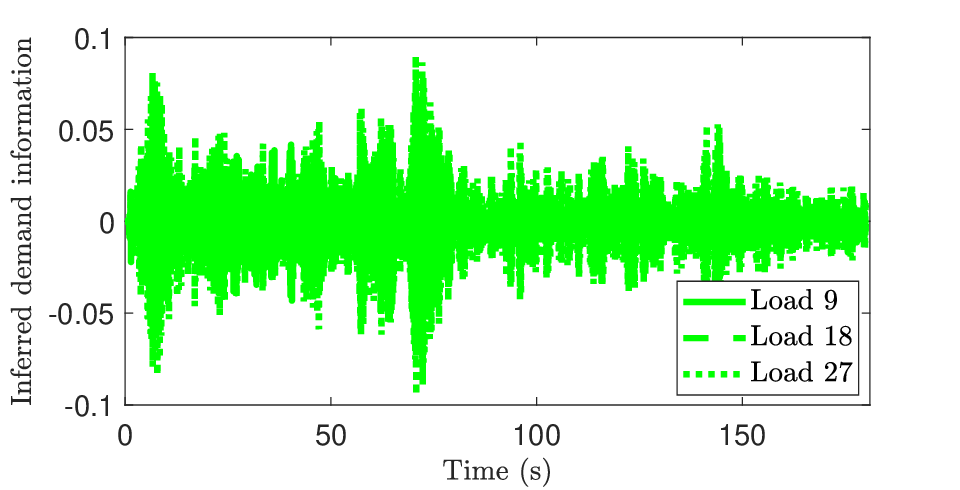}
\vspace{-4mm}
\caption{Inferred demand information on loads $9, 18$ and $27$ at bus $2$ using the power command trajectories for the following two control schemes:   (top) Extended Primal-Dual scheme, and (bottom) Privacy-Preserving scheme.}
\label{Privacy_B}
\vspace{-5mm}
\end{figure}

To demonstrate that the \textit{Privacy-Preserving scheme} ensures the privacy of the prosumption profiles against intelligent eavesdroppers, we considered an observer scheme that aims to infer the controllable demand  using a model of the power command dynamics and knowledge of the power command signals.
In particular, by evaluating the power command derivative and the value of $\psi$, an eavesdropper may attempt to observe the generation and controllable demand profiles by reversing \eqref{5b}, i.e. using $ \widetilde{s} =  \Gamma \dot{p}^c + H \psi$. 
Figure \ref{Privacy_B} demonstrates the result from such observer scheme for the same three loads considered in Fig. \ref{Privacy_A}, when the \textit{Extended Primal-Dual} and  \textit{Privacy-Preserving} schemes are implemented.
From Fig. \ref{Privacy_B}, it follows that an intelligent eavesdropper may obtain the   controllable demand profiles when the \textit{Extended Primal-Dual scheme} is applied. By contrast, the application of the \textit{Privacy-Preserving scheme} ensures that the demand is private against intelligent eavesdroppers, since the retrieved information is distorted by the signal $n_{k,j}$.
{The outcomes of Section \ref{V} are synopsized in Table \ref{tab_synopsis}, which provides a comparison between the properties of each of the considered schemes in terms of enabling convergence, an optimal power allocation and prosumption privacy against naive and intelligent eavesdroppers.}

 \begin{table}[h]
 \centering
 \normalsize
\begin{tabular}{|l|l|l|l|l|}
 \hline
 Property /  Scheme & (i)   
 & (ii)  
 & (iii)
  & (iv) 
   \\
 \hline
Convergence & $\checkmark$ & $\checkmark$  & $\checkmark$ & $\checkmark$ \\
 \hline
Optimality & \xmark & $\checkmark$ & $\checkmark$ & $\checkmark$ \\
 \hline
Privacy -  & $\checkmark$  & \xmark & $\checkmark$  & $\checkmark$ \\
 Naive Eavesdroppers &   &  &  & \\
 \hline
Privacy -  & $\checkmark$  &   \xmark & \xmark & $\checkmark$ \\
 Intelligent Eavesdroppers  &  &  &  & \\
 \hline
\end{tabular}
\caption{{Comparison of the properties of: (i) the Integral Action scheme, (ii) the \textit{Primal-Dual} scheme, (iii) the \textit{Extended Primal-Dual} scheme, and (iv) the \textit{Privacy-Preserving} scheme, following the outcomes presented in Section \ref{V}.}}
\vspace{-3mm}
 \label{tab_synopsis}
\end{table}

{
\subsection{Statistical validation}\label{sec_stats}

This section aims to offer statistical validation to our simulation results.
In our simulations we randomly selected
(a)~the number of controllable loads at each of 20 load buses, 
(b)~the loads where the step change in demand was applied,
(c)~the cost coefficient parameters $q_{k,j}, k \in \mathcal{N}_j, j \in \mathcal{N}$ and the values for $\xi_{k,j}(t)$ and $n^f_{k,j}(t)$ associated with scheme (iv) and
(d)~the structure of the (connected) communication graph for schemes (iii) and (iv).

To validate our results for a broad range of randomly selected parameters, we repeated the presented simulations  $100$ times for each scheme and compared their performance on the following three properties, (a) the frequency at the end of the simulation duration, (b) the synchronization of the marginal cost values at the end of the simulation duration and (c) the ability of eavesdroppers to infer prosumption profiles. 
The results are shown in Table \ref{tab_statistics}.
Rows 1 and 2 of Table \ref{tab_statistics} demonstrate that in all cases the frequency converges to its nominal value, although  schemes (iii) and (iv) result in slower frequency convergence.
Nevertheless, in all cases the frequencies converge to within $0.015$ Hz from the nominal value within three minutes, which is well within the secondary frequency control timeframe.
To produce the third and fourth rows of Table \ref{tab_statistics}, we removed the top and bottom $5\%$ of the marginal cost values, to avoid  outliers.
We then obtained the variance and maximum difference from the average value of the marginal cost values based on this set.
The results demonstrate that scheme (i) yields significantly larger values of these two quantities
 compared to schemes (ii), (iii) and (iv).
The latter validates that scheme (i) yields a suboptimal response, while schemes (ii), (iii) and (iv) yield an optimal response, as follows from the very low variance and maximum difference-from-mean marginal cost values.
Lastly, to validate the results in Fig. \ref{Privacy_A} and \ref{Privacy_B}, we considered the ability of eavesdroppers to infer the prosumption profile values.
It should be noted that this aim is not relevant for scheme (i) since it involves no communication.
In addition, for scheme (ii), since prosumption values are communicated, the error is always zero.
For schemes (iii) and (iv), the approach to infer the prosumption profiles was by reversing \eqref{5b}, i.e. using $ \widetilde{s} =  \Gamma \dot{p}^c + H \psi$, and calculating the mean absolute prosumption inference error.
The results demonstrate that for scheme (iii) this approach always allows to infer prosumption while for scheme (iv) there exists a significant mean inference error.
Hence, the presented statistical results validate the theoretical results  in this paper.

 \begin{table}[ht!]
 \centering
 \normalsize
\begin{tabular}{|l|c|c|c|c|}
 \hline
 Property /  Scheme & (i) 
 & (ii)  
 & (iii) 
  & (iv) 
   \\
 \hline
Mean freq. (Hz) & \hspace{-1mm}-0.0001 & \hspace{-1mm}-0.0001  & \hspace{-1mm}-0.012 & \hspace{-1mm}-0.013 \hspace{-3mm}\\
 \hline
Min freq. (Hz) & \hspace{-1mm}-0.001 & \hspace{-1mm}-0.0003 & \hspace{-1mm}-0.014 & \hspace{-1mm}-0.015 \hspace{-3mm}\\
 \hline
Marginal Cost -  &  \multirow{2}{*}{\hspace{-1mm}2$\cdot 10^{-5}$}& \multirow{2}{*}{\hspace{-1mm}2$\cdot 10^{-8}$} & \multirow{2}{*}{$10^{-6}$}  & \multirow{2}{*}{\hspace{-1mm}7$\cdot 10^{-7}$ \hspace{-3mm}}
\\
Variance  &  &  &  & \\
 \hline
Marginal Cost -    &  &  &  & \\
Max Difference & \hspace{-1mm}0.012  &   \hspace{-1mm}3$\cdot 10^{-4}$& \hspace{-1mm}2$\cdot 10^{-3}$ & \hspace{-1mm}2$\cdot 10^{-3}$ \hspace{-3mm}\\
from Mean  &  &  &  & \\
 \hline
 Mean Absolute   &  &  &  & \\ 
 Prosumption &  N/A & 0 & 0 & {\hspace{-1mm}0.14 \hspace{-3mm}} \\
   Inference Error  &  &  &  & \\
 \hline
\end{tabular}
\caption{{Statistical comparison from simulations implementing each of the following  control schemes $100$ times: (i)~Integral action scheme, (ii)~\textit{Primal-Dual} scheme, (iii)~\textit{Extended Primal-Dual} scheme, and (iv)~\textit{Privacy-Preserving} scheme.
The compared quantities are the average and minimum frequency at the end of the simulation time duration, the marginal cost variance, the maximum difference between the average and any non-outlier marginal cost values, and the mean absolute prosumption inference error.}}
\vspace{-3mm}
 \label{tab_statistics}
\end{table}
}

\section{Conclusion}\label{VI}
We have considered the problem of  enabling an optimal power allocation and simultaneously preserving the privacy of generation and controllable demand profiles within the secondary frequency control timeframe.
To enhance the intuition on our results,
two types of eavesdroppers were defined; naive eavesdroppers that {do not possess/make use of knowledge of the internal {system} dynamics to} analyze the intercepted signals and intelligent eavesdroppers that
use knowledge of the underlying dynamics to infer the privacy-sensitive prosumption profiles.
We proposed the \textit{Extended Primal-Dual scheme}, which implements a controller at each privacy-seeking unit in the power grid to provide improved privacy properties.
The proposed scheme enables privacy guarantees against naive eavesdroppers.
However,  the generation/demand profiles may be inferred  by intelligent eavesdroppers   using the communicated signal trajectories and information on the underlying {system} dynamics.
To resolve this issue,  we proposed the \textit{Privacy-Preserving scheme}, which shares the structure of the \textit{Extended Primal-Dual scheme} but also incorporates a privacy-enhancing signal at each controller.
The latter {continuously} adjusts the  response speed  of the controllers, making model based observations inaccurate, and  disturbs the generation/demand profile information within the controllers, enabling privacy against intelligent eavesdroppers.
For both proposed schemes, we  provide analytic stability, optimality and privacy guarantees.
Our presented results are distributed, locally verifiable and applicable to  general network configurations.
 The applicability of the proposed schemes is demonstrated with simulations on the NPCC 140-bus system where we  show that stability is preserved, and improved privacy properties and an optimal power allocation are attained.

\balance

\vspace{-2mm}
\bibliography{bib_file}

\end{document}